\documentclass[
 aip,
 cha,
 amsmath,
 amssymb,
 reprint,
 longbibliography
]{revtex4-2}

\usepackage{graphicx}
\usepackage{dcolumn}
\usepackage{bm}
\usepackage{amsmath}
\usepackage{amssymb}
\usepackage{xcolor}
\usepackage{dcolumn}
\usepackage{bm}

\usepackage[utf8]{inputenc}
\usepackage[T1]{fontenc}
\usepackage{mathptmx}
\usepackage{etoolbox}

\usepackage{tabularx}
\usepackage{array}
\usepackage{placeins}

\usepackage{amsthm}
\usepackage{color}
\usepackage{comment}
\usepackage{enumitem}

\usepackage{booktabs}

\usepackage{hyperref}
\usepackage{multirow}

\newcommand{\myEdit}[1]{{\color{black} #1}}

\begin{document}

\title{Collective Phase Reorganization and Cluster Synchronization in Networks of Coupled Gumowski--Mira Maps}

\author{Hammed O. Fatoyinbo}
\email{hammed.fatoyinbo@aut.ac.nz}
\affiliation{Department of Mathematical Sciences, School of Engineering, Computer and Mathematical Sciences, Auckland University of Technology, Auckland 1010, New Zealand}

\author{Indranil Ghosh}
\email{indra.ghosh@ucd.ie}
\affiliation{School of Mathematics and Statistics, University College Dublin, Dublin D04 V1W8, Ireland}

\date{\today}

\begin{abstract}
We investigate the collective dynamics of networks composed of diffusively coupled Gumowski–Mira maps and determine how changes in the intrinsic dynamics of the local map influence the organization of network states. The coupling strength and the local control parameter are treated as bifurcation parameters, and the collective dynamics are characterized using the largest Lyapunov exponent, a synchronization error, cluster number, and a combined state classification. Two local parameter regimes are examined under the same ring topology and coupling rule. In both cases, the network supports synchronized fixed points, periodic cluster states, chaotic cluster states, and chaotic-incoherent dynamics, although the locations and distributions of these states in the collective parameter plane differ. Long-run simulations confirm all four representative collective states in both local parameter regimes. Complete Lyapunov-spectrum calculations show that the synchronized fixed points and periodic clusters are nonchaotic, whereas the chaotic clusters are hyperchaotic, with 19 and 21 positive Lyapunov exponents in Cases~A and~B, respectively. The chaotic-incoherent states are also hyperchaotic and have 50 positive exponents in both cases. A multiple-initial-condition analysis further shows that the robustness of a collective classification depends on the selected parameter region. These results demonstrate that variations in the intrinsic map dynamics modify both the spatial organization and the high-dimensional instability of collective states in coupled discrete-time networks. 
\end{abstract}


\maketitle

\begin{quote}


\textcolor{black}{Networks of interacting nonlinear units arise across science and engineering, from neural and ecological systems to power grids and coordinated technological networks. A fundamental challenge is to understand how changes in the behaviour of individual components influence the large-scale organization of the entire system. This work shows that collective behaviour is not determined by network structure and interaction strength alone: changes in the intrinsic dynamics of the individual units can substantially reorganize synchronization, clustering, incoherence, and instability, even when the network architecture remains unchanged. The results provide a general framework for distinguishing temporal instability from spatial organization and demonstrate that highly unstable dynamics can coexist with coherent group structure. These insights contribute to a broader understanding of how microscopic dynamics shape macroscopic patterns in complex nonlinear networks.}
\end{quote}

\section{Introduction} \label{sec:I}

The dynamical properties of complex systems represent a central theme in the study of nonlinear phenomena, with applications spanning physics, engineering, biology, ecology, neuroscience, and the social sciences~\cite{Strogatz2000,Boccaletti2006,Arenas2008,Fatoyinbo22}. A particularly significant subclass of such systems is formed by networks whose nodes evolve according to nonlinear dynamical laws and whose interactions give rise to emergent collective behavior over time. Depending on the local node dynamics, the coupling topology, and the interaction strength, networks of nonlinear oscillators can display a broad spectrum of macroscopic states, including complete synchronization, cluster synchronization, incoherent dynamics, multistability, chimera-like patterns, and collective chaos~\cite{Pecora1998,Kuramoto2002,Abrams2004,Majhi2019,Bick2020}.

Coupled map lattices and coupled map networks constitute a versatile and powerful framework for the analysis of such phenomena. In contrast to continuous-time dynamical systems, discrete-time maps are computationally efficient and thereby enable extensive explorations of high-dimensional parameter spaces. Consequently, they are particularly well suited for identifying stability regions, bifurcation structures, synchronization domains, spatiotemporal pattern formation, and transitions between regular and irregular regimes of dynamics~\cite{Kaneko1989,Kaneko1992,Bunimovich1988}. Coupled maps have furthermore been extensively employed to investigate complete synchronization, partial synchronization, cluster synchronization, and the structure and stability of invariant synchronization manifolds~\cite{Belykh2000,Zhuang2002,Baraviera2023}. These characteristics render coupled map networks a natural and mathematically tractable framework for studying how local nonlinear dynamics and network-mediated interactions jointly determine emergent collective behavior.

Discrete-time map-based models provide a useful framework for studying collective dynamics in nonlinear networks because they combine computational simplicity with rich intrinsic dynamics. Networks of coupled maps can exhibit synchronization transitions, multistability, coherent--incoherent coexistence, cluster formation, and irregular or chaotic collective behavior~\cite{Ibarz2011, Muni2022,Ghosh2023,Nair2024}. Such systems therefore provide a computationally tractable setting for investigating how interactions between nonlinear units generate macroscopic dynamical organization. Previous studies of networks composed of map-based oscillators have shown that bifurcation structures, Lyapunov-based chaos diagnostics, complexity measures, and synchronization indices can be jointly employed to characterize collective dynamics in discrete-time nonlinear networks~\cite{Nair2024,Ghosh2025}.

In the present work, we choose the Gumowski--Mira map as the local node dynamics because it is a two-dimensional discrete-time system whose intrinsic dynamics vary strongly and nonmonotonically with its control parameters ~\cite{GumowskiMira1980,Mira1996,GumowskiMira2006}. Within a single mathematical model, parameter variation can produce fixed points, periodic orbits, invariant-curve dynamics, intermittent behavior, and chaotic attractors. This pronounced parameter-dependent dynamical reorganization is the principal feature that makes the map suitable for the present study. It allows us to compare how qualitatively different local dynamical regimes reorganize collective network states while the network topology, coupling rule, parameter domain, and numerical protocol are kept unchanged. Previous studies have examined intermittency and synchronization in coupled Gumowski--Mira maps and have characterized its attractors using Lyapunov-based indicators ~\cite{Ambika2006,Saha2006}. However, the present work addresses a different question: how changes in the intrinsic regime of the local map alter the organization of synchronized, clustered, and incoherent states in the collective parameter plane.

Synchronization in networks of nonlinear oscillators has been extensively investigated, as it constitutes a fundamental mechanism by which interacting units coordinate and organize their dynamics~\cite{Pecora1998,Arenas2008}. In the most elementary scenario, all oscillators asymptotically converge to an identical trajectory, yielding complete synchronization. In many networked systems, however, global synchronization fails to occur. Instead, the nodes may segregate into subsets that synchronize internally, giving rise to cluster synchronization~\cite{Belykh2000,Bick2020,Pecora2014,Sorrentino2016}. Cluster synchronization is therefore a central form of partial organization in complex networks, where coherent groups coexist while remaining dynamically distinct from one another. Clustered states may also occur in irregular or chaotic regimes, leading to chaotic cluster synchronization in coupled-map and networked-map systems~\cite{Bohle2021}. In other parameter regimes, coherent and incoherent subpopulations may coexist, leading to chimera-like or partially synchronized states~\cite{Kuramoto2002,Abrams2004,Omelchenko2015,Hart2016}. These partial collective states are particularly pertinent in networks of coupled maps, because low-dimensional maps can generate complex temporal dynamics while remaining computationally tractable for large-scale parameter explorations.

\textcolor{black}{Discrete-time neuron models provide one important application of coupled-map networks because they combine computational efficiency with nonlinear excitable dynamics~\cite{Chialvo1995,Kanagaraj2024,Used2024}. In the present work, the Gumowski–Mira map is treated as an abstract two-dimensional nonlinear oscillator for investigating how parameter-dependent changes in local dynamics influence network-level synchronization, clustering, and incoherence.}

The present study aligns with this broader line of research but addresses a distinct problem. Rather than solely determining whether a coupled Gumowski--Mira network achieves synchronization, we investigate how the intrinsic dynamics of the local oscillators restructure the global collective phase space. This distinction is crucial, since synchronization thresholds alone do not provide a complete characterization of the emergent behavior of nonlinear networks. For example, a system may exhibit periodic yet incoherent dynamics, irregular or chaotic clustered organization, or a highly structured partition into multiple stable groups. To systematically capture these possibilities, we employ a combination of dynamical and collective observables: the largest Lyapunov exponent to distinguish regular and irregular regimes, the synchronization error to quantify global coherence, and a cluster-count measure to estimate the number of coherent groups within the network~\cite{Eckmann1985,Politi2013,Vadivasova2016}.

We investigate a ring network of diffusively coupled Gumowski--Mira maps and construct two-parameter \myEdit{collective-state} diagrams in the plane spanned by the local control parameter \myEdit{$\mu$} and coupling strength \myEdit{$\epsilon$}. \myEdit{Two intrinsic parameter regimes of the isolated map are considered. The resulting network diagrams reveal synchronized, clustered, and incoherent states whose locations and distributions depend on the local parameters $(a,b)$. Representative long-time simulations confirm synchronized fixed points, periodic two-cluster states, chaotic two-cluster states, and chaotic-incoherent dynamics in both local parameter regimes. Complete Lyapunov-spectrum calculations further show that the synchronized fixed points and periodic clusters are nonchaotic, whereas the chaotic clusters and chaotic-incoherent states are hyperchaotic.}

\myEdit{ In this work, \textit{collective phase organization} refers to the arrangement and distribution of synchronized, clustered, and incoherent collective states in the $(\mu,\epsilon)$ parameter plane. The term \texttt{phase} therefore denotes a numerically classified collective dynamical regime rather than the instantaneous phase angle of an individual oscillator.}

The primary contribution of this work is \myEdit{a comparative characterization of the} collective-\myEdit{state} organization \myEdit{generated by two distinct intrinsic regimes of the Gumowski--Mira map. By keeping the topology, coupling rule, parameter domain, and numerical protocol fixed, we determine how changes in the local map parameters modify the organization and robustness of synchronized, clustered, and incoherent network states. We additionally use the complete Lyapunov spectrum to distinguish nonchaotic cluster synchronization from hyperchaotic cluster synchronization and to quantify the higher-dimensional instability of chaotic-incoherent dynamics.} 

The remainder of the paper is organized as follows. In Sec.~\ref{sec:II}, we introduce the coupled Gumowski--Mira network and \myEdit{characterize the two local dynamical regimes.} In Sec.~\ref{sec:III}, we \myEdit{define the largest Lyapunov exponent, synchronization error, cluster number, and collective-state classification. Sections~\ref{sec:IV} and~\ref{sec:V} present the collective parameter spaces for Cases~A and~B, respectively}. In Sec.~\ref{sec:VI}, we compare the two regimes and discuss the impact of local oscillator dynamics on the emergent network behavior. \textcolor{black}{Section~\ref{sec:representative_collective_states} presents representative collective states, examines their complete Lyapunov spectra, and evaluates sensitivity to initial conditions, finite observation intervals, numerical thresholds, and network size.} Finally, Sec.~\ref{sec:VIII} summarizes the main results and outlines potential directions for future research.

\section{Coupled Gumowski--Mira Network} \label{sec:II}

We consider a network of $N$ identical Gumowski--Mira maps arranged on a ring with nearest-neighbor coupling. The isolated Gumowski--Mira map is given by
\begin{align}
x_{n+1} &= ay_n(1-by_n^2)+y_n+F(x_n,\mu), \label{eq:gm_x}\\
y_{n+1} &= -x_n+F(x_{n+1},\mu), \label{eq:gm_y}
\end{align}
where
\begin{equation}
F(x,\mu)=\mu x+\frac{2x\textcolor{blue}{^2}(1-\mu)}{1+x^2}.
\label{eq:F}
\end{equation}
Here, $x_n$ and $y_n$ are the state variables at discrete time $n$, while $a$, $b$, and $\mu$ are control parameters. The parameter $\mu$ is used as one of the principal bifurcation parameters in this work.

To study collective behavior, we couple $N$ Gumowski--Mira maps through the $x$ variable. Let $(x_i(n),y_i(n))$ denote the state of node $i$ at time $n$, where $i=1,2,\ldots,N$. Periodic boundary conditions are imposed so that
\begin{equation}
x_{0}(n)=x_N(n), \qquad x_{N+1}(n)=x_1(n).
\end{equation}
For each node, we first define the local Gumowski--Mira update
\begin{equation}
u_i(n)=ay_i(n)\left[1-by_i^2(n)\right]+y_i(n)+F(x_i(n),\mu).
\label{eq:local_update}
\end{equation}
The coupled network is then written as
\begin{align}
x_i(n+1) &=
(1-\epsilon)u_i(n)
+\frac{\epsilon}{2}
\left[x_{i-1}(n)+x_{i+1}(n)\right],
\label{eq:coupled_x}\\
y_i(n+1) &=
-x_i(n)+F(x_i(n+1),\mu),
\label{eq:coupled_y}
\end{align}
where $\epsilon$ is the coupling strength. The case $\epsilon=0$ corresponds to uncoupled Gumowski--Mira maps, while increasing $\epsilon$ strengthens the influence of nearest neighbors. The form of Eq.~\eqref{eq:coupled_x} implements a diffusive nearest-neighbor interaction through a convex combination of the local update and neighboring states.

\textcolor{black}{The nearest-neighbour ring was selected as a canonical local-coupling architecture widely used to investigate spatiotemporal dynamics, pattern formation, and synchronization in coupled-map systems~\cite{Kaneko1992,Kaneko1993,KanekoYanagita2014}. Periodic boundary conditions eliminate edge effects, while the regular ring provides a controlled baseline for isolating the combined influence of the intrinsic Gumowski–Mira dynamics and local diffusive coupling. More general coupling ranges and network architectures may introduce additional effects that are beyond the scope of the present study \cite{KoillerYoung2010, Atay2004}.}



\subsection{Selection of local dynamical regimes}

The two local parameter regimes considered in this work are selected to represent distinct intrinsic dynamics of the isolated Gumowski--Mira map. This choice is motivated by the strong parameter sensitivity of the Gumowski--Mira system, which is known to exhibit periodic, quasiperiodic, intermittent, and chaotic behavior as its control parameters are varied~\cite{GumowskiMira1980,Mira1996,GumowskiMira2006,Ambika2006}. By comparing two different local regimes under the same network topology and coupling rule, we can isolate the influence of the local oscillator dynamics on the emergent collective phase organization.

We focus on two representative regimes:
\begin{align}
\text{Case A:} \qquad & a=0.01, \qquad b=0.05 \label{eq:caseA} \\
\text{Case B:} \qquad & a=0.33, \qquad b=0.10 \label{eq:caseB}
\end{align}
For both cases, the local control parameter is varied in the interval
\begin{equation}
\mu\in[-0.6,0.6],
\end{equation}
while the coupling strength is varied in
\begin{equation}
\epsilon\in[0,0.30].
\end{equation}
Thus, the explored collective parameter plane is
\begin{equation}
\label{eq:pplane}
(\mu,\epsilon)\in[-0.6,0.6]\times[0,0.30].
\end{equation}



 \subsection{\myEdit{Dynamics of the isolated Gumowski--Mira map}} \label{sec:isolated_dynamics} 
 \myEdit{ Before examining the coupled network, we characterize the dynamics of the isolated Gumowski--Mira map over the complete interval $\mu\in[-0.6,0.6]$ used in the collective parameter-space analysis. This characterization is necessary because the principal objective of the present study is to determine how changes in the intrinsic node dynamics influence the organization of collective states. In particular, a bifurcation diagram alone may not reliably distinguish a densely sampled invariant curve from a chaotic attractor. We therefore supplement the bifurcation diagrams with the largest Lyapunov exponent, $\lambda_{\max}^{\mathrm{iso}}$, of the isolated two-dimensional map. For the nonlinear function defined in Eq.~(3), the derivative with respect to $x$ is 
 
 \begin{equation} F_x(x,\mu) = \mu+ \frac{4(1-\mu)x}{(1+x^2)^2}. 
 \label{eq:isolated_derivative}
 \end{equation} 
 Defining 
 \begin{equation*} A_n=F_x(x_n,\mu), \qquad B_n=1+a\left(1-3b y_n^2\right), \qquad C_n=F_x(x_{n+1},\mu), 
 \end{equation*} 
 
 the one-step Jacobian matrix of the isolated map is 
 \begin{equation} 
 \mathbf{J}_n = \begin{pmatrix} A_n & B_n \\ -1+C_nA_n & C_nB_n \end{pmatrix}. 
 \label{eq:isolated_jacobian} 
 \end{equation} 
 
 The largest Lyapunov exponent $\lambda_{\max}^{\mathrm{iso}}$ is estimated by propagating a tangent vector under $\mathbf{J}_n$, normalizing the tangent vector after each iteration, and averaging the corresponding logarithmic growth factors after discarding the transient dynamics. We classify the isolated dynamics using the numerical tolerance $\lambda_{\mathrm{tol}}=10^{-3}$. Specifically, parameter values satisfying 
 
 \begin{equation*} \lambda_{\max}^{\mathrm{iso}}<-\lambda_{\mathrm{tol}}, \qquad \left|\lambda_{\max}^{\mathrm{iso}}\right| \leq\lambda_{\mathrm{tol}}, \qquad \lambda_{\max}^{\mathrm{iso}}>\lambda_{\mathrm{tol}} 
 \end{equation*} 
 
 are classified as attracting regular, numerically neutral, and chaotic, respectively. A numerically neutral largest Lyapunov exponent may be associated with invariant-curve dynamics, weakly attracting periodic behavior, or proximity to a local bifurcation. Consequently, numerically neutral values are interpreted together with the corresponding phase portraits and are not, by themselves, regarded as sufficient evidence of quasiperiodicity. For the isolated-map calculations, the interval $\mu\in[-0.6,0.6]$ was sampled using 1201 uniformly spaced parameter values, corresponding to $\Delta\mu=10^{-3}$. At each value of $\mu$, the map was initialized independently from the fixed initial condition $(x_0,y_0)=(0,0.5)$. The first $n_0^{\mathrm{iso}}=5000$ iterations were discarded to remove transient behavior. The subsequent 300 iterates were retained for constructing the bifurcation diagrams, while the largest Lyapunov exponent was estimated over $T_{\mathrm{iso}}=5000$ iterations. Each value of $\mu$ was initialized from the same fixed initial condition rather than using parameter continuation. } 

 \myEdit{ The isolated-map characteristics of the two parameter cases are summarized in Table~\ref{tab:local_regimes}. Because the network topology, coupling rule, explored parameter plane, and numerical protocol are kept fixed, varying the local parameters $(a,b)$ provides a controlled means of examining how different intrinsic dynamical organizations are reflected in the coupled network. } 

\begin{table*}[t!]
\centering
\caption{Local parameter regimes used to compare the collective dynamics.} 
\label{tab:local_regimes} \begin{tabular}{@{} c c c l @{}} 
\hline 
Regime & $a$ & $b$ & \multicolumn{1}{c}{Isolated-map characteristic} \\ 
\hline 
Case~A & $0.01$ & $0.05$ & \parbox[t]{8.2cm}{Alternating chaotic and attracting regular intervals, with narrow numerically neutral regions. } \\[1.5ex]
Case~B & $0.33$ & $0.10$ & \parbox[t]{8.2cm}{ Distinct chaotic and strongly attracting regular intervals, with invariant-curve candidates. } \\ 
\hline 
\end{tabular} 
\end{table*}
 
 \begin{figure}[t!] 
 \centering
 \includegraphics[width=\columnwidth]{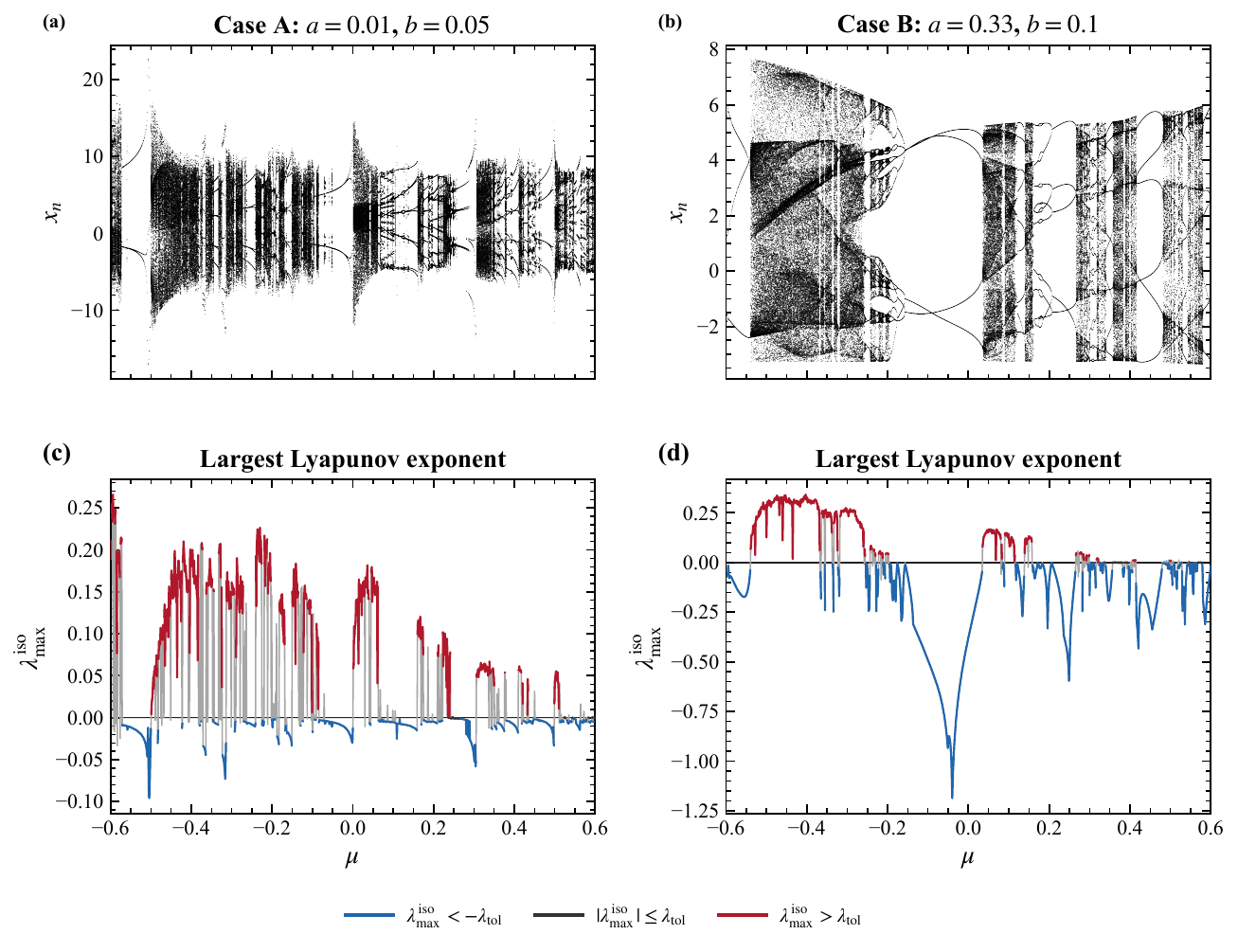} 
 \caption{ Dynamics of the isolated Gumowski--Mira map over the parameter interval used in the network analysis. Panels~(a) and~(b) show the bifurcation diagrams for Cases~A and~B, respectively, while panels~(c) and~(d) show the corresponding largest Lyapunov exponents $\lambda_{\max}^{\mathrm{iso}}$. Blue, gray, and red portions denote $\lambda_{\max}^{\mathrm{iso}}<-\lambda_{\mathrm{tol}}$, $|\lambda_{\max}^{\mathrm{iso}}| \leq\lambda_{\mathrm{tol}}$, and $\lambda_{\max}^{\mathrm{iso}}>\lambda_{\mathrm{tol}}$, respectively, with $\lambda_{\mathrm{tol}}=10^{-3}$. The gray horizontal band marks the numerically neutral interval. For each value of $\mu$, the first $n_0^{\mathrm{iso}}=5000$ iterations were discarded, the next 300 iterates were used for the bifurcation diagram, and the exponent was estimated over $T_{\mathrm{iso}}=5000$ iterations. } 
 \label{fig:local_regimes} 
 \end{figure}

  \myEdit{ Figure~\ref{fig:local_regimes} presents the bifurcation diagrams and largest Lyapunov exponents for Cases~A and~B. Case~A exhibits repeated alternation among positive-exponent chaotic intervals, attracting regular windows, and comparatively narrow numerically neutral regions. Case~B displays a different ordering and distribution of chaotic, strongly attracting regular, and numerically neutral intervals, including positive-$\mu$ values for which the orbit approaches a smooth closed curve. Thus, the two cases differ in the organization of their intrinsic dynamics throughout the investigated $\mu$ interval, rather than only in the geometry of selected attractors. 
  
  Previous analysis of the Gumowski--Mira map has identified saddle--node, flip, and Neimark--Sacker bifurcations, together with periodic stability regions and tongue-like structures ~\cite{BorinOliveira2026}. These results provide a useful analytical context for interpreting changes in the isolated-map dynamics. However, several windows in Fig.~\ref{fig:local_regimes} involve higher-period attractors, and we do not assign a specific bifurcation mechanism to every observed transition. The regimes considered here are therefore identified using the combined bifurcation-diagram, Lyapunov-exponent, and phase-portrait evidence. 
 
} 
 
 \begin{figure}[t!] 
 \centering \includegraphics[width=\columnwidth]{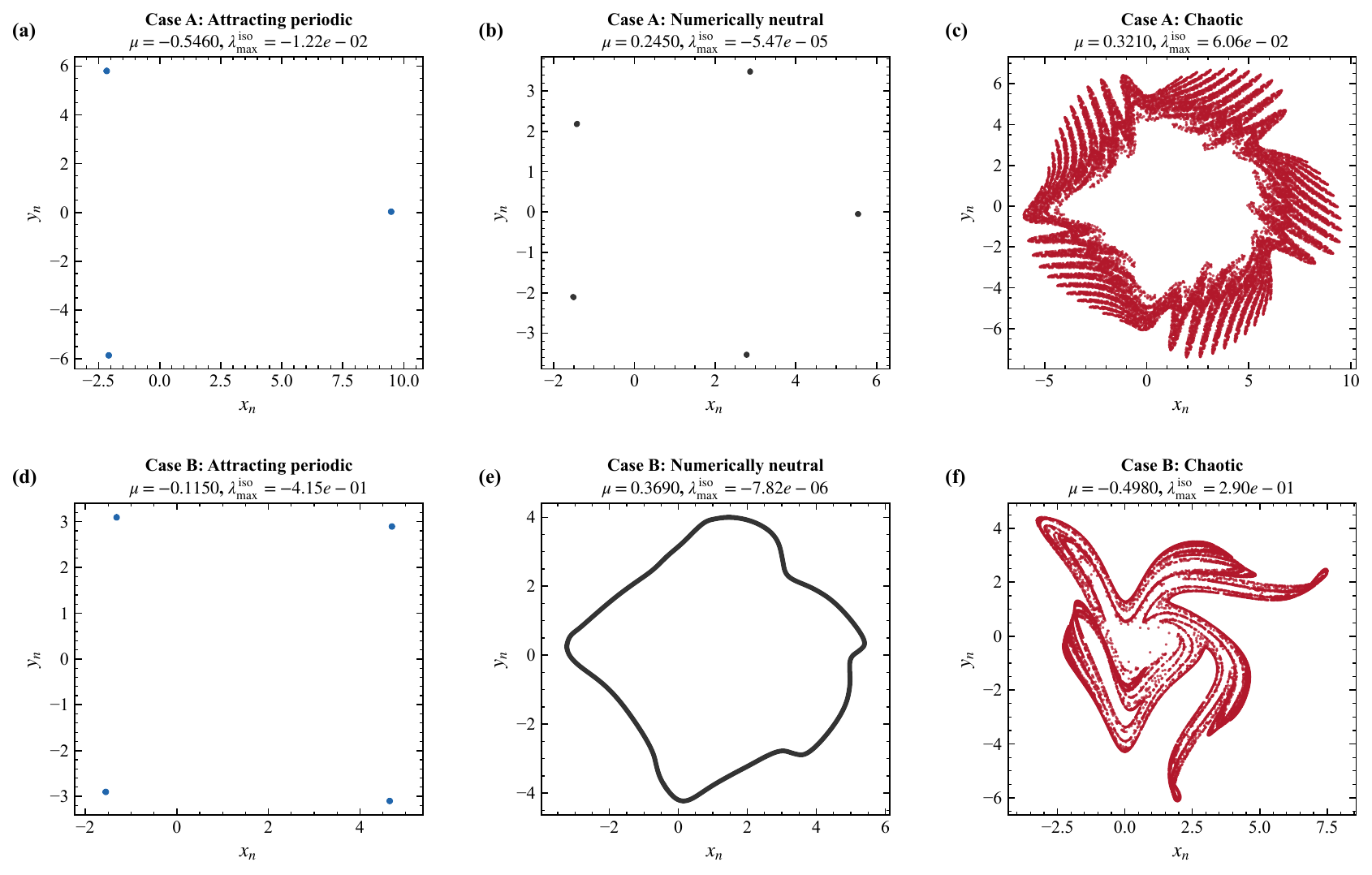} 
 \caption{Representative phase portraits of the isolated Gumowski--Mira map selected from the bifurcation and largest-Lyapunov-exponent analyses. Panels~(a)--(c) show attracting periodic, numerically neutral, and chaotic states for Case~A, respectively. Panels~(d)--(f) show the corresponding dynamical classes for Case~B. The parameter value and largest Lyapunov exponent are reported above each panel. The smooth closed curve in panel~(e) is consistent with invariant-curve dynamics. Numerically neutral classification alone is not regarded as sufficient evidence of quasiperiodicity. Each representative state was recomputed using $5\times10^{4}$ transient iterations and $5\times10^{4}$ iterations for the final Lyapunov estimate. The phase portraits were constructed from $10^{4}$ post-transient points. } 
 \label{fig:phase_portraits} 
 \end{figure} 
 
 \myEdit{ For each local parameter case, the illustrated values of $\mu$ in Fig.~\ref{fig:phase_portraits} were selected from the combined bifurcation and Lyapunov analyses to represent attracting periodic, numerically neutral, and chaotic dynamics. The numerically neutral orbit for Case~B forms a smooth closed curve, which is consistent with invariant-curve dynamics. By contrast, the numerically neutral example for Case~A is interpreted more cautiously because a near-zero largest Lyapunov exponent does not, by itself, establish quasiperiodicity. These isolated-map diagnostics provide the local dynamical basis for the subsequent comparison of the two coupled networks. }



\section{Numerics}\label{sec:III}

To characterize the collective dynamics of the coupled Gumowski--Mira network, we use four complementary diagnostics: the largest Lyapunov exponent $\lambda_{\max}$, the synchronization error $E$, the detected cluster number \(N_c\), and a combined collective-state classification. This follows the general numerical philosophy used in previous studies, where multiple indicators are combined to distinguish regularity, chaos, coherence, synchronization, and cluster formation in coupled nonlinear systems~\cite{Eckmann1985,Politi2013,Nair2024,Ghosh2025}. Detailed definitions of the numerical diagnostics and implementation procedures are provided in Appendix~\ref{app:appendix}.

\myEdit{Unless otherwise stated, the network contains $N=50$ nodes. The value \(N=50\) is treated as the baseline finite network size. To assess finite-size dependence, the representative parameter pairs were additionally simulated for \(N=25\) and \(N=100\), using ten independently generated initial network states for each parameter pair and network size. The collective parameter plane defined in Eq.~\eqref{eq:pplane} was discretized using a $200\times200$ uniform grid. For each parameter pair, the node states were initialized independently according to 
\begin{equation} x_i(0),y_i(0)\sim\mathcal{U}[-1,1], \qquad i=1,\ldots,N. 
\label{eq:network_initial_conditions} 
\end{equation} 
The first $n_0=1500$ iterations were discarded. The synchronization error and cluster statistics were then evaluated over $T=2000$ iterations, while the largest Lyapunov exponent was estimated over $T_{\lambda}=3000$ post-transient iterations. At each parameter pair, all diagnostics were calculated from the same initial network state and reference trajectory.

The numerical thresholds used throughout the network analysis are 
\[ \lambda_{\mathrm{tol}}=10^{-3},\qquad E_{\mathrm{tol}}=10^{-8},\qquad \eta=5\times10^{-2},\qquad N_c^*=3. \] 

The cluster number was evaluated at 20 uniformly distributed post-transient times and summarized by the rounded median. Consequently, the cluster number reported for a representative state describes its typical post-transient organization and need not equal the instantaneous number of amplitude groups visible in a single final network snapshot.

At each sampled iteration, the node values \(x_i(n)\) are sorted in ascending order, and a new cluster is introduced whenever the gap between two adjacent sorted values exceeds the threshold \(\eta=0.05\). No additional minimum cluster size is imposed by this threshold procedure; therefore, isolated nodes may occur as singleton clusters in the raw cluster count. However, a state is classified as cluster synchronization only when
\[ 2\leq N_c\leq N_c^*. \]
Thus, a state with \(N_c=1\) cannot receive a periodic-, neutral-, or chaotic-cluster classification.

The largest Lyapunov exponent \(\lambda_{\max}\) is used to distinguish regular and irregular temporal dynamics. Negative values of $\lambda_{\max}$ indicate asymptotically regular dynamics, whereas positive values indicate sensitive dependence on initial conditions.}

\textcolor{black}{The sign of \(\lambda_{\max}\) characterizes the temporal dynamics but does not, by itself, establish synchronization.
Complete synchronization is confirmed separately by
\[
E\leq E_{\mathrm{tol}} \qquad\text{and}\qquad N_c=1.
\]
Thus, \(\lambda_{\max}\), \(E\), and \(N_c\) provide complementary
temporal and spatial diagnostics.} Since the computed exponents are finite-time estimates, values close to zero are treated as marginal and are not overinterpreted. We therefore introduce a small tolerance $\lambda_{\rm tol}$ when assigning collective states.


The synchronization error $E$ measures the dispersion of node states around the network mean and is used to quantify global coherence. Complete synchronization corresponds to $E$ approaching zero. The cluster number $N_c$ estimates the number of coherent groups in the network and is used to distinguish complete synchronization from cluster synchronization and highly fragmented states.

The collective-state classification used in the phase diagrams is summarized in Appendix Table~\ref{tab:collective_classification}. Complete synchronization is identified by a single coherent group and vanishing synchronization error. Periodic cluster synchronization corresponds to regular dynamics in which the network separates into multiple coherent groups. Chaotic cluster synchronization denotes clustered collective dynamics with positive finite-time largest Lyapunov exponent. \myEdit{States satisfying $N_c>N_c^*=3$ are classified as incoherent, or highly fragmented, irrespective of the sign of the finite-time largest Lyapunov exponent. An incoherent state with $\lambda_{\max}>\lambda_{\mathrm{tol}}$ is described more specifically as chaotic incoherent, whereas one with $\lambda_{\max}<-\lambda_{\mathrm{tol}}$ is regular but highly fragmented.}

The classification should be interpreted as a finite-time numerical classification based on the largest Lyapunov exponent, synchronization error, and cluster-count statistics. In particular, the label ``chaotic cluster synchronization'' denotes a clustered state with positive finite-time largest Lyapunov exponent. Since finite-time estimates may be sensitive near transition boundaries, selected representative states are further examined using the complete Lyapunov spectrum in Sec.~\ref{sec:representative_collective_states}.

For this supplementary check, \myEdit{the complete Lyapunov spectrum was computed for each of the eight representative states shown in Figs.~\ref{fig:representative_states_caseA} and \ref{fig:representative_states_caseB}. \myEdit{The complete spectrum was obtained by evolving a complete tangent-space basis together with the reference trajectory and applying QR orthonormalization after each iteration~\cite{Benettin1980a, Benettin1980b}.} For $N=50$, the tangent space has dimension $2N=100$, and each calculation therefore yields 100 Lyapunov exponents. The reference trajectory was evolved for 5000 transient iterations, after which the spectrum was accumulated over 20000 iterations. Intermediate estimates were recorded after 2000, 5000, 10000, 15000, and 20000 accumulation iterations to monitor convergence.} 

Let $n_+$ denote the number of Lyapunov exponents satisfying
\begin{equation*}
\lambda_i>10^{-3}.
\end{equation*}
A representative state is classified as \myEdit{nonchaotic} when $n_+=0$, chaotic when $n_+=1$, and hyperchaotic when $n_+\geq2$. The criterion $n_+\geq2$ follows the standard definition of hyperchaos as dynamics with at least two positive Lyapunov exponents~\cite{Letellier2007}. \myEdit{ The complete spectrum provides the long-run dynamical classification of each representative trajectory, while the synchronization error and rounded median cluster number describe its spatial organization.}

\textcolor{black}{To assess the robustness of the finite-time classification, additional calculations were performed using multiple independent initial network states and discrete choices of the post-transient observation interval, Lyapunov tolerance, synchronization tolerance, and cluster threshold. The tested values and corresponding results are reported in Sec.~\ref{sec:representative_collective_states}.}

\section{Case A: Collective Phase Organization}\label{sec:IV}

We first investigate the collective dynamics of the coupled Gumowski--Mira network for the local parameter regime
\begin{equation*}
a=0.01,\qquad b=0.05.
\end{equation*}
This regime is referred to as Case A. The parameter $\mu$ is varied in the interval $[-0.6,0.6]$, while the coupling strength $\epsilon$ is varied in the interval $[0,0.30]$. For each parameter pair $(\mu,\epsilon)$, we compute the largest Lyapunov exponent, the synchronization error, the number of clusters, and the corresponding collective-state classification.

\subsection{Collective Lyapunov parameter space}

\begin{figure}[h] 
\centering 
\includegraphics[width=\columnwidth]{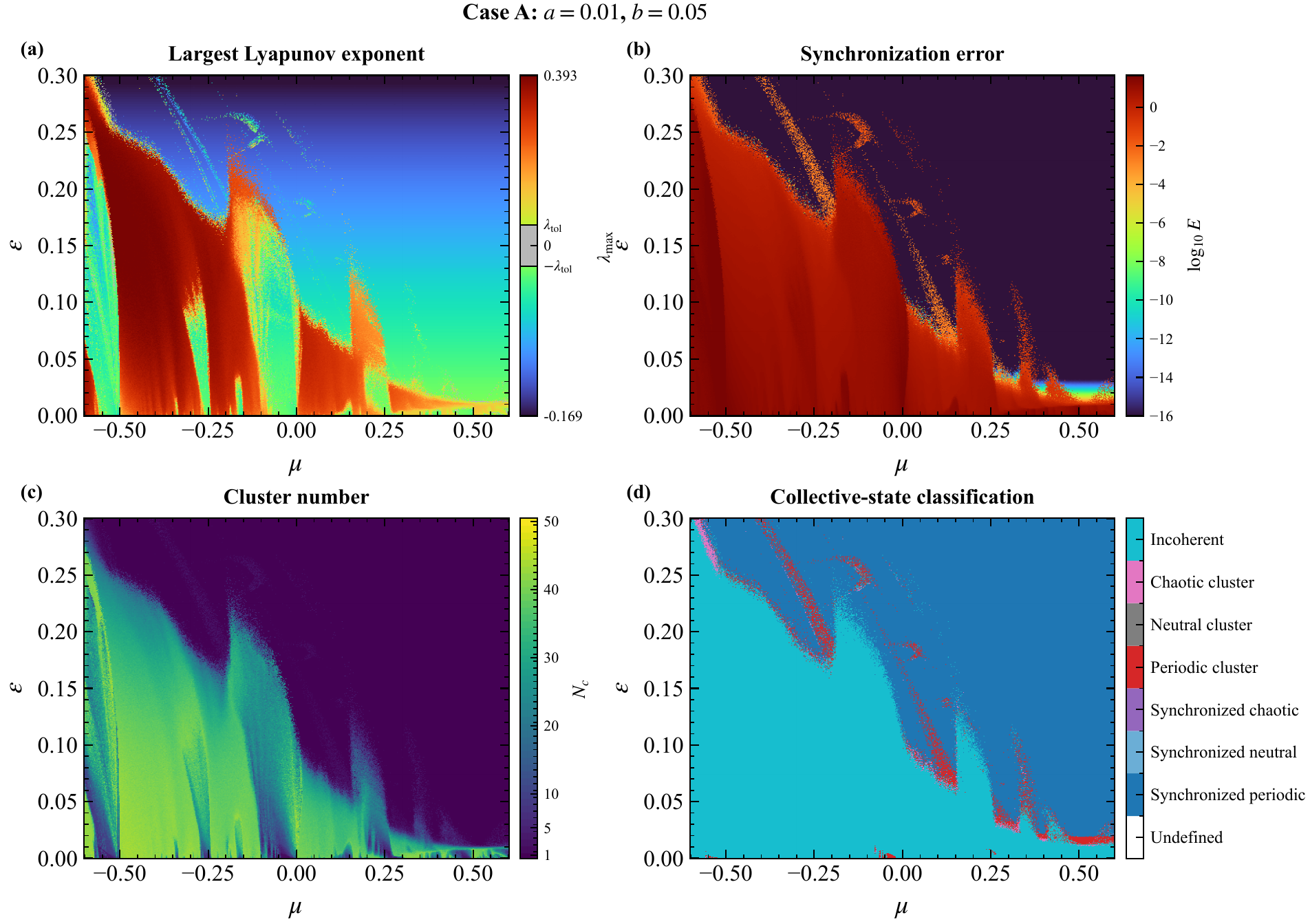} 
\caption{ Collective parameter-space organization of the coupled Gumowski--Mira network for Case~A, with $a=0.01$ and $b=0.05$. Panels~(a)--(d) show the largest Lyapunov exponent, time-averaged synchronization error, rounded median cluster number, and combined collective-state classification, respectively, in the $(\mu,\epsilon)$ plane. Gray denotes $|\lambda_{\max}|\leq\lambda_{\mathrm{tol}}$ in panel~(a), while panel~(c) uses discrete integer color intervals for $N_c$, with $N_c=1$ denoting complete synchronization. The network contains $N=50$ nodes, and the parameter plane is sampled on a $200\times200$ uniform grid. The largest Lyapunov exponent, synchronization error, and cluster number at each parameter pair were calculated from the same initial network state and reference trajectory. Numerically neutral dynamics satisfy $|\lambda_{\max}|\leq10^{-3}$. } 
\label{fig:case_a_parameter_spaces} 
\end{figure}

Figure~\ref{fig:case_a_parameter_spaces}(a) shows the largest Lyapunov exponent in the $(\mu,\epsilon)$ parameter plane for Case~A. \myEdit{Positive Lyapunov exponents occupy substantial portions of the weak-coupling region, whereas negative exponents become more prevalent as the coupling strength increases. The transition is nonuniform across the investigated $\mu$ interval.}

\myEdit{For more negative values of $\mu$, the upper envelope of the positive-exponent region exhibits an overall downward shift as $\mu$ increases. For $\mu\gtrsim-0.15$, the organization becomes more fragmented, with alternating bands and localized pockets of positive, neutral, and negative Lyapunov exponents. The diagram therefore exhibits an overall sloping transition region rather than a single smooth boundary extending across the complete parameter interval. Nevertheless, the increasing prevalence of negative Lyapunov exponents at stronger coupling indicates that coupling generally regularizes the collective dynamics in Case~A. }

\subsection{Synchronization parameter space}

To determine whether the regular regions identified by the Lyapunov exponent correspond to globally coherent dynamics, we compute the time-averaged synchronization error $E$. Figure~\ref{fig:case_a_parameter_spaces}(b) shows the synchronization error in the same $(\mu,\epsilon)$ parameter plane. 

\myEdit{ The synchronization-error diagram exhibits an overall transition toward lower values of $E$ as the coupling strength increases. The location and structure of this transition depend on $\mu$. For more negative values of $\mu$, the boundary of the low-error region broadly follows the sloping organization observed in the Lyapunov-exponent diagram. For $\mu\gtrsim-0.15$, the transition becomes more fragmented and includes narrow bands and localized pockets of reduced synchronization error. We therefore describe this feature as a nonuniform synchronization domain rather than a single smooth synchronization boundary. } 

The synchronization error does not vanish throughout the complete regular region identified in Fig.~\ref{fig:case_a_parameter_spaces}(a). Regular collective dynamics therefore do not necessarily imply complete synchronization. Instead, the network may organize into several coherent groups, motivating the cluster-number analysis below.

\subsection{Cluster synchronization parameter space}

Figure~\ref{fig:case_a_parameter_spaces}(c) shows the rounded median cluster number $N_c$ in the $(\mu,\epsilon)$ parameter plane. \myEdit{This diagnostic distinguishes complete synchronization from low-cluster organization and highly fragmented dynamics. 

The cluster-number distribution broadly follows the organization observed in the synchronization-error diagram. Larger cluster numbers occur predominantly at weaker coupling, whereas stronger coupling generally reduces the number of detected clusters. The transition to low-cluster organization is nonuniform and contains narrow bands and localized parameter-dependent structures, particularly for $\mu\gtrsim-0.15$. 

The low-cluster region includes both complete synchronization, corresponding to $N_c=1$, and states composed of a small number of distinct coherent groups. The network therefore does not undergo a simple transition directly from incoherence to complete synchronization. Cluster synchronization instead constitutes an important and extended form of collective organization.}


\subsection{Collective phase diagram}

 To jointly account for the temporal and spatial properties, we classify each point in the $(\mu,\epsilon)$ parameter plane using the criteria summarized in Appendix Table~\ref{tab:collective_classification}. The resulting collective-state diagram is presented in Fig.~\ref{fig:case_a_parameter_spaces}(d). 
 
 \myEdit{ The diagram confirms an overall transition from weak-coupling incoherent and clustered dynamics toward more coherent collective states at stronger coupling. This transition is nonuniform and is interrupted by narrow periodic, neutral, and chaotic-cluster structures. These localized structures correspond to the bands and pockets visible in the Lyapunov-exponent, synchronization-error, and cluster-number panels. } 
 
 Periodic cluster synchronization occupies substantial portions of the parameter plane, particularly at intermediate coupling strengths. Synchronized-periodic dynamics become prominent in parts of the stronger-coupling region. Regions classified as chaotic cluster synchronization occur mainly within localized transition structures where the finite-time largest Lyapunov exponent is positive while the rounded median cluster number remains small. Incoherent states occur predominantly at weaker coupling and within fragmented regions near collective-state boundaries. 
 
 \myEdit{ The higher-resolution parameter scan resolves narrow chaotic-cluster structures within and near the transition region. These states satisfy $\lambda_{\max}>\lambda_{\mathrm{tol}}$, $E>E_{\mathrm{tol}}$, and $2\leq N_c\leq N_c^*$, demonstrating that irregular temporal dynamics can coexist with a small number of coherent groups. The complete Lyapunov spectrum of the selected Case~A chaotic-cluster representative is examined in Sec.~\ref{sec:representative_collective_states}. } 
 
 In summary, Case~A exhibits a \myEdit{broadly organized but nonuniform collective parameter space. Increasing coupling generally promotes temporal regularity, reduces the synchronization error, and decreases the rounded median cluster number. However, the transition proceeds through parameter-dependent bands and localized clustered states rather than across one smooth boundary. The resolved chaotic-cluster regions additionally show that coupling-induced spatial organization can persist even when the collective temporal dynamics remain irregular. }

\section{Case B: Collective Phase Organization}\label{sec:V}

\myEdit{We next consider Case~B, with $a=0.33$ and $b=0.10$. The network topology, coupling rule, parameter domain, initialization protocol, and numerical settings are identical to those used for Case~A. Consequently, differences between the two collective parameter spaces arise from changing the intrinsic local parameters $(a,b)$.}

\subsection{Collective Lyapunov parameter space}

\begin{figure}[h] 
\centering 
\includegraphics[width=\columnwidth]{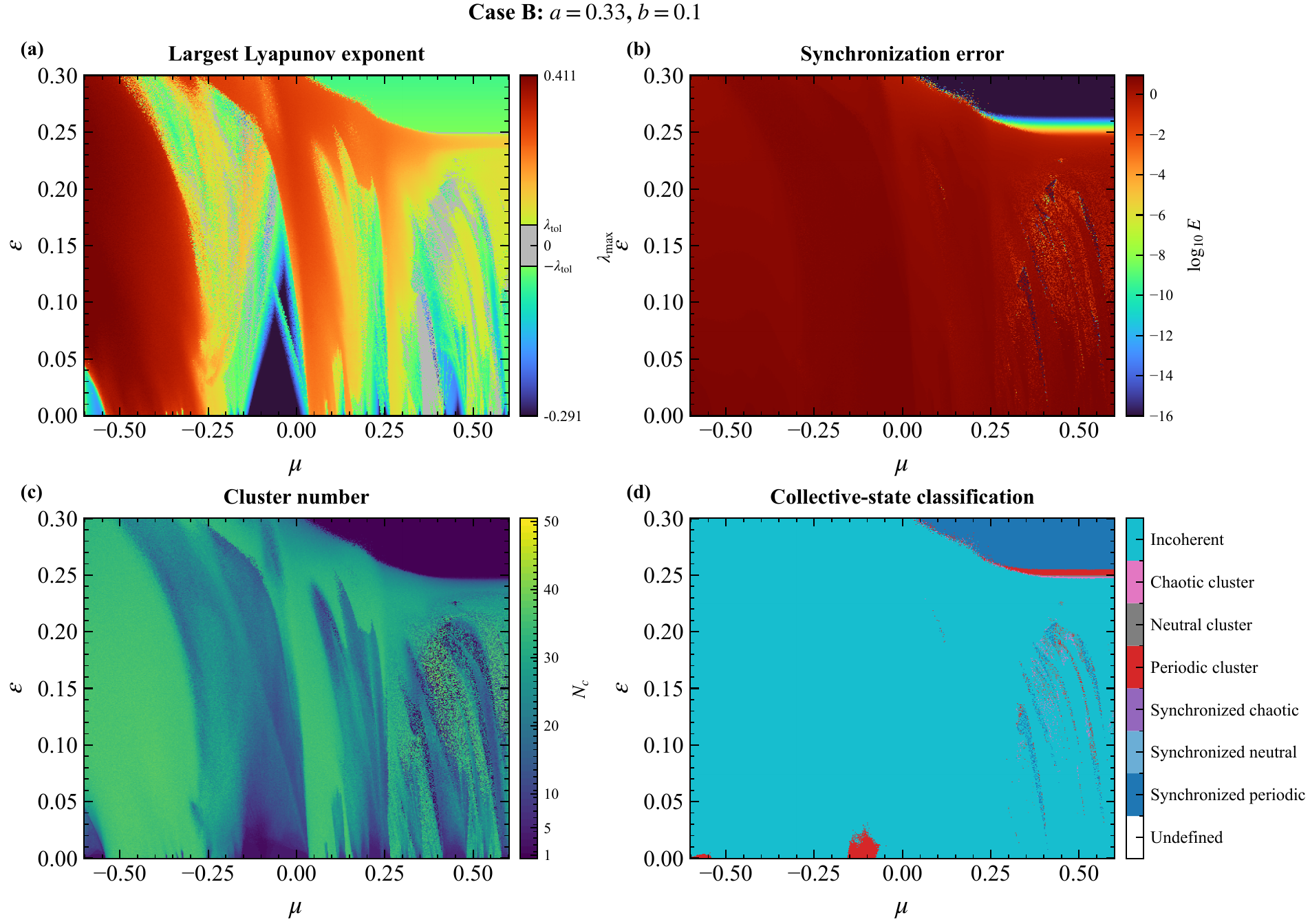} 
\caption{ Collective parameter-space organization of the coupled Gumowski--Mira network for Case~B, with $a=0.33$ and $b=0.10$. Panels~(a)--(d) show the largest Lyapunov exponent, time-averaged synchronization error, rounded median cluster number, and combined collective-state classification, respectively, in the $(\mu,\epsilon)$ plane. Gray denotes $|\lambda_{\max}|\leq\lambda_{\mathrm{tol}}$ in panel~(a), while panel~(c) uses discrete integer color intervals for $N_c$, with $N_c=1$ denoting complete synchronization. The network contains $N=50$ nodes, and the parameter plane is sampled on a $200\times200$ uniform grid. The largest Lyapunov exponent, synchronization error, and cluster number at each parameter pair were calculated from the same initial network state and reference trajectory. Numerically neutral dynamics satisfy $|\lambda_{\max}|\leq10^{-3}$. } 
\label{fig:case_b_parameter_spaces} 
\end{figure}

Figure~\ref{fig:case_b_parameter_spaces}(a) shows the largest Lyapunov exponent in the $(\mu,\epsilon)$ parameter plane for Case~B. \myEdit{Positive finite-time Lyapunov exponents occupy extensive portions of the parameter plane, indicating that irregular collective dynamics persist over a broad range of coupling strengths.

The distribution is strongly nonuniform. Regions with negative or numerically neutral largest Lyapunov exponents form curved bands, narrow tongues, and localized pockets within the broader positive-exponent background. A pronounced regular region occurs near $\mu\approx-0.1$ at weak and intermediate coupling, while additional alternating regular, neutral, and irregular structures occur for positive $\mu$. At stronger coupling and sufficiently positive $\mu$, the largest Lyapunov exponent becomes negative over an extended domain. Thus, Case~B does not exhibit a single transition boundary. Instead, the parameter plane contains an interwoven arrangement of regular, neutral, and irregular collective regimes. }

\subsection{Synchronization parameter space}

Figure~\ref{fig:case_b_parameter_spaces}(b) shows the time-averaged synchronization error for Case~B. \myEdit{Much of the parameter plane has a comparatively large synchronization error, indicating that the network generally remains outside the complete-synchronization manifold. 

 A prominent low-error domain appears at strong coupling and positive $\mu$, particularly for approximately $\mu\gtrsim0.2$ and $\epsilon\gtrsim0.25$. Its lower boundary is curved and parameter dependent rather than uniform. Additional narrow bands and localized low-error structures occur at smaller coupling strengths, especially within the positive-$\mu$ region. The synchronization organization is therefore more heterogeneous than in Case~A and consists of one extended strong-coupling domain together with finer localized structures. } 

The synchronization map also shows that temporal regularity does not necessarily imply complete synchronization. Several regions with negative or near-zero largest Lyapunov exponent retain a nonvanishing synchronization error. These regions may therefore correspond to regular clustered states rather than globally synchronized motion.

\subsection{Cluster synchronization parameter space}

\myEdit{Figure~\ref{fig:case_b_parameter_spaces}(c) shows the rounded median cluster number $N_c$. Large cluster numbers occur over substantial portions of the parameter plane, consistent with the comparatively large synchronization error observed in Fig.~\ref{fig:case_b_parameter_spaces}(b). 

 Low-cluster regions occur most prominently at strong coupling and positive $\mu$, where they overlap with the extended low-error domain. Additional low-cluster bands and pockets occur at weaker coupling, including structures near $\mu\approx-0.1$ and within the positive-$\mu$ range. The transition between high- and low-cluster organization is irregular and contains several narrow parameter-dependent structures.

The low-cluster regions include both complete synchronization, corresponding to $N_c=1$, and states composed of a small number of coherent groups. Consequently, the reduction of the synchronization error and cluster number does not represent a single direct transition from incoherence to complete synchronization. Instead, the network passes through a heterogeneous collection of synchronized, clustered, and fragmented states.}

\subsection{Collective phase diagram}

Using the classification criteria summarized in Appendix Table~\ref{tab:collective_classification}, we construct the collective-state diagram shown in Fig.~\ref{fig:case_b_parameter_spaces}(d). 

\myEdit{Incoherent dynamics occupy a substantial portion of the parameter plane, particularly at weak and intermediate coupling. An extended synchronized-periodic domain occurs at strong coupling and positive $\mu$, in agreement with the negative largest Lyapunov exponent, low synchronization error, and small cluster number observed in panels~(a)--(c). Its boundary is curved and contains narrow periodic-cluster and chaotic-cluster structures. } 

Periodic-cluster states occur in localized regions, including a regular domain near $\mu\approx-0.1$ at weak coupling and narrow structures near collective-state boundaries. Chaotic-cluster states also occur in localized bands and pockets where the finite-time largest Lyapunov exponent is positive but the rounded median cluster number remains small. 

\myEdit{The parameter scan shows a representative chaotic-cluster state at $(\mu,\epsilon)=(0.5351,0.2477)$. Although its finite-time largest Lyapunov exponent is comparatively small and positive, the complete spectrum reported in Sec.~\ref{sec:representative_collective_states} confirms that the selected trajectory possesses multiple positive Lyapunov exponents. } 

\myEdit{ Overall, Case~B exhibits a heterogeneous collective parameter space dominated by incoherent dynamics over broad regions, with an extended synchronized-periodic domain at strong coupling and positive $\mu$. Regular, neutral, periodic-cluster, and chaotic-cluster states form curved bands, narrow tongues, and localized pockets within and near the boundaries of these broader domains. Compared with Case~A, the collective organization is more fragmented and more strongly dependent on both $\mu$ and $\epsilon$. The selected chaotic-cluster state further demonstrates that low-cluster spatial organization can coexist with high-dimensional temporal instability. }

\section{Comparison of the Two Local Regimes}\label{sec:VI}

The results in Secs.~IV and~V show that the same coupled network architecture can generate qualitatively different collective phase organizations depending on the intrinsic dynamics of the local Gumowski--Mira oscillator. In both cases, the network consists of identical nearest-neighbor coupled maps arranged on a ring, and the same parameter plane $(\mu,\epsilon)$ is explored. Therefore, the differences between Cases~A and~B arise from the change in the local parameters $(a,b)$. 

\subsection{Comparison of collective transition structures} 

In Case~A, with 
\begin{equation}
a=0.01,\qquad b=0.05,
\end{equation} 
\myEdit{the parameter plane exhibits an overall sloping but nonuniform transition toward regular and coherent collective dynamics as the coupling strength increases. This broad organization is visible in the Lyapunov-exponent, synchronization-error, cluster-number, and collective-state diagrams, although the transition contains narrow bands and localized parameter-dependent structures.}

By contrast, Case~B, with
\begin{equation}
a=0.33,\qquad b=0.10, 
\end{equation} 
\myEdit{exhibits a more heterogeneous collective organization. Broad incoherent regions coexist with curved regular and neutral bands, localized cluster states, and an extended synchronized-periodic domain at strong coupling and positive $\mu$. }

This comparison demonstrates that \myEdit{ changes in the intrinsic local dynamics do not simply shift a single synchronization threshold. Instead, they modify the geometry and distribution of the collective states throughout the $(\mu,\epsilon)$ plane. The two cases nevertheless support the same four representative categories: synchronized fixed points, periodic clusters, chaotic clusters, and chaotic-incoherent states.}

\subsection{Robustness of periodic cluster synchronization} 

Despite the strong differences between the two cases, one feature remains robust: the prevalence of periodic cluster synchronization. In both Case~A and Case~B, substantial portions of the parameter plane are classified as periodic clusters. These states correspond to regular temporal dynamics together with multiple coherent groups, namely
\begin{equation} 
\lambda_{\max}<-\lambda_{\mathrm{tol}}, \qquad E>E_{\mathrm{tol}}, \qquad 2\leq N_c\leq N_c^*. 
\label{eq:periodic_cluster_condition} 
\end{equation} 

This observation shows that complete synchronization is not the generic outcome of coupling in the Gumowski--Mira network. Instead, the network frequently self-organizes into stable cluster states. Such behavior is consistent with coupling suppressing temporal irregularity without necessarily forcing all nodes onto one globally synchronized trajectory. 

The persistence of periodic cluster synchronization across both local regimes suggests that clustering is a robust emergent property of the coupled Gumowski--Mira network. The details of the cluster organization, however, depend strongly on the intrinsic dynamics of the local map. \myEdit{The selected periodic-cluster representatives illustrate this distinction. Both representatives have rounded median cluster number $N_c=2$, and all 100 computed Lyapunov exponents are negative, but their instantaneous spatial arrangements differ substantially. Thus, the periodic-cluster classification identifies a common form of regular partial synchronization without requiring identical cluster geometry in the two local parameter regimes.} 

\subsection{Hyperchaotic clustered states and incoherent domains} 

Another important difference between the two cases concerns the distribution of irregular clustered states and incoherent dynamics. In Case~A, regions classified as chaotic cluster synchronization are comparatively localized and occur mainly within transition structures separating incoherent dynamics, periodic cluster synchronization, and complete synchronization. In Case~B, chaotic-cluster states occur in localized bands and pockets whose positions and geometries differ from those in Case~A. These states satisfy 
\begin{equation} 
\lambda_{\max}>\lambda_{\mathrm{tol}}, \qquad E>E_{\mathrm{tol}}, \qquad 2\leq N_c\leq N_c^*, 
\label{eq:chaotic_cluster_condition} 
\end{equation} 
indicating that the network retains a small number of coherent groups while the finite-time largest Lyapunov exponent is positive. 

\myEdit{ The complete Lyapunov spectra confirm that the selected chaotic-cluster states are hyperchaotic. The Case~A and Case~B representatives have 19 and 21 exponents satisfying $\lambda_i>10^{-3}$, respectively. Their leading exponents are $\lambda_1=1.592\times10^{-1}$ for Case~A and $\lambda_1=2.604\times10^{-3}$ for Case~B. Although the leading exponent for Case~B is comparatively close to the adopted threshold, the presence of 21 positive exponents confirms that its hyperchaotic classification is not based on a single marginally unstable direction. } 

Incoherent states occupy broader regions in Case~B, particularly at weak and intermediate coupling and near collective-state boundaries. Under the adopted classification, a highly fragmented state satisfies 
\begin{equation}
N_c>N_c^*.
\label{eq:incoherent_condition} 
\end{equation} 

An incoherent state is described more specifically as chaotic incoherent when it also satisfies 
\begin{equation} 
\lambda_{\max}>\lambda_{\mathrm{tol}}. 
\label{eq:chaotic_incoherent_condition} 
\end{equation} 

\myEdit{The selected chaotic-incoherent representatives have 50 positive Lyapunov exponents in both cases. Their instability is therefore distributed across more tangent-space directions than in the chaotic-cluster representatives. This provides a dynamical complement to the spatial classification: chaotic clusters retain a small rounded median cluster number, whereas chaotic-incoherent states are highly fragmented. } 

\myEdit{Overall, the same coupling architecture supports regular, synchronized, clustered, hyperchaotic, and incoherent states in both local parameter regimes, but their locations and geometries depend strongly on the intrinsic map parameters. Case~A exhibits an overall sloping but nonuniform transition toward regular and coherent dynamics, whereas Case~B displays a more heterogeneous organization with broader incoherent regions and an extended synchronized-periodic domain at strong coupling and positive $\mu$. These results demonstrate that spatial organization and temporal instability are related but distinct properties: a network may retain low-cluster organization while remaining unstable along many independent tangent-space directions. }

\section{Representative Collective States} \label{sec:representative_collective_states} 
\begin{figure}[!t] 
\centering 
\includegraphics[width=\columnwidth]{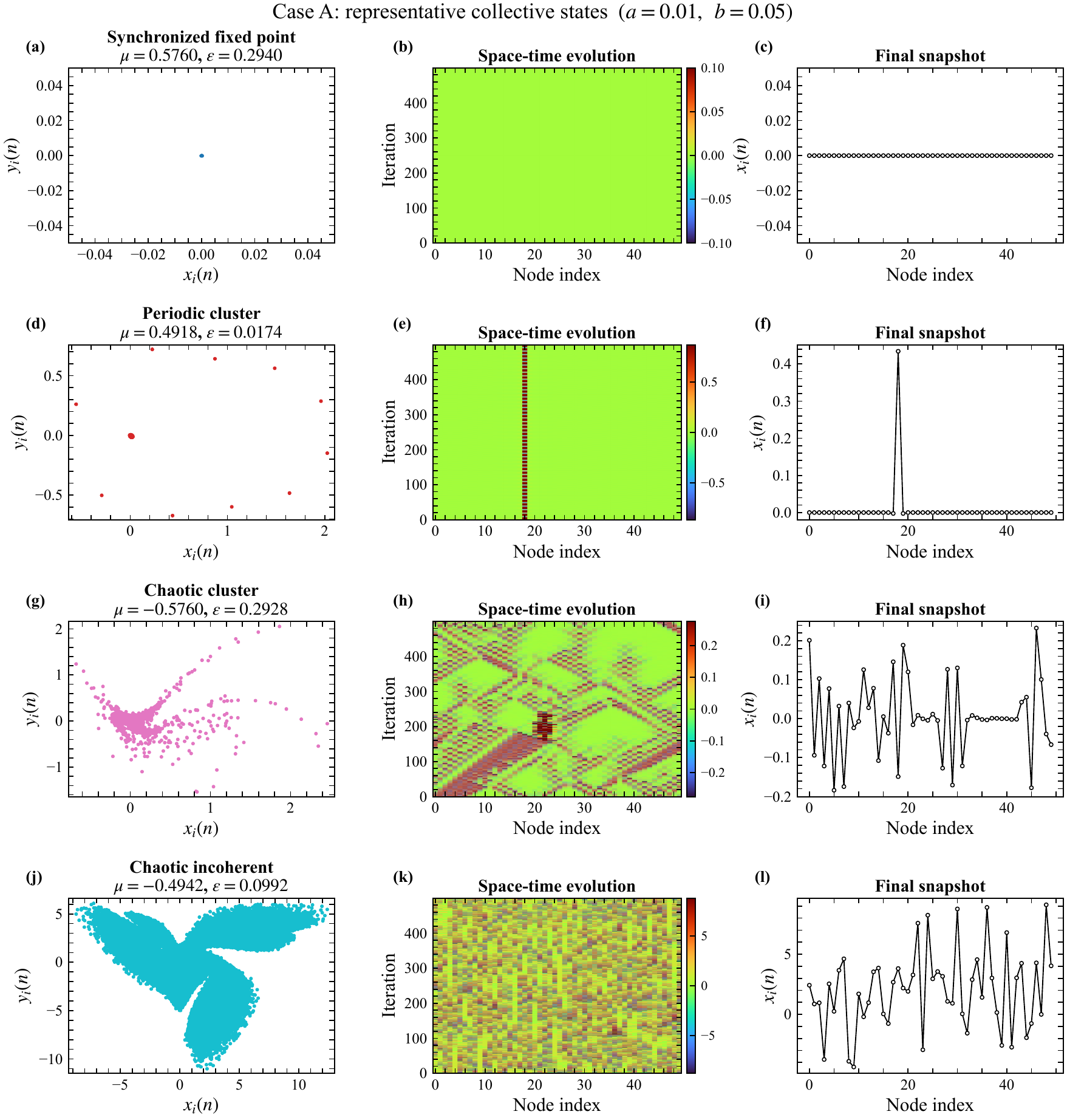} 
\caption{Representative collective states of the coupled Gumowski--Mira network for Case~A, with $a=0.01$ and $b=0.05$. \myEdit{ The rows show a synchronized fixed point at $(\mu,\epsilon)=(0.5760,0.2940)$, a periodic cluster state at $(\mu,\epsilon)=(0.4918,0.0174)$, a chaotic cluster state at $(\mu,\epsilon)=(-0.5760,0.2928)$, and a chaotic-incoherent state at $(\mu,\epsilon)=(-0.4942,0.0992)$. The first column shows the physical $(x_i,y_i)$ phase-space trajectories of all nodes, the second column shows the final 500 iterations of the space-time evolution, with color representing the instantaneous value of $x_i(n)$}, and the third column shows the final network snapshot.}
\label{fig:representative_states_caseA} 
\end{figure}

\begin{figure}[!t] 
\centering 
\includegraphics[width=\columnwidth]{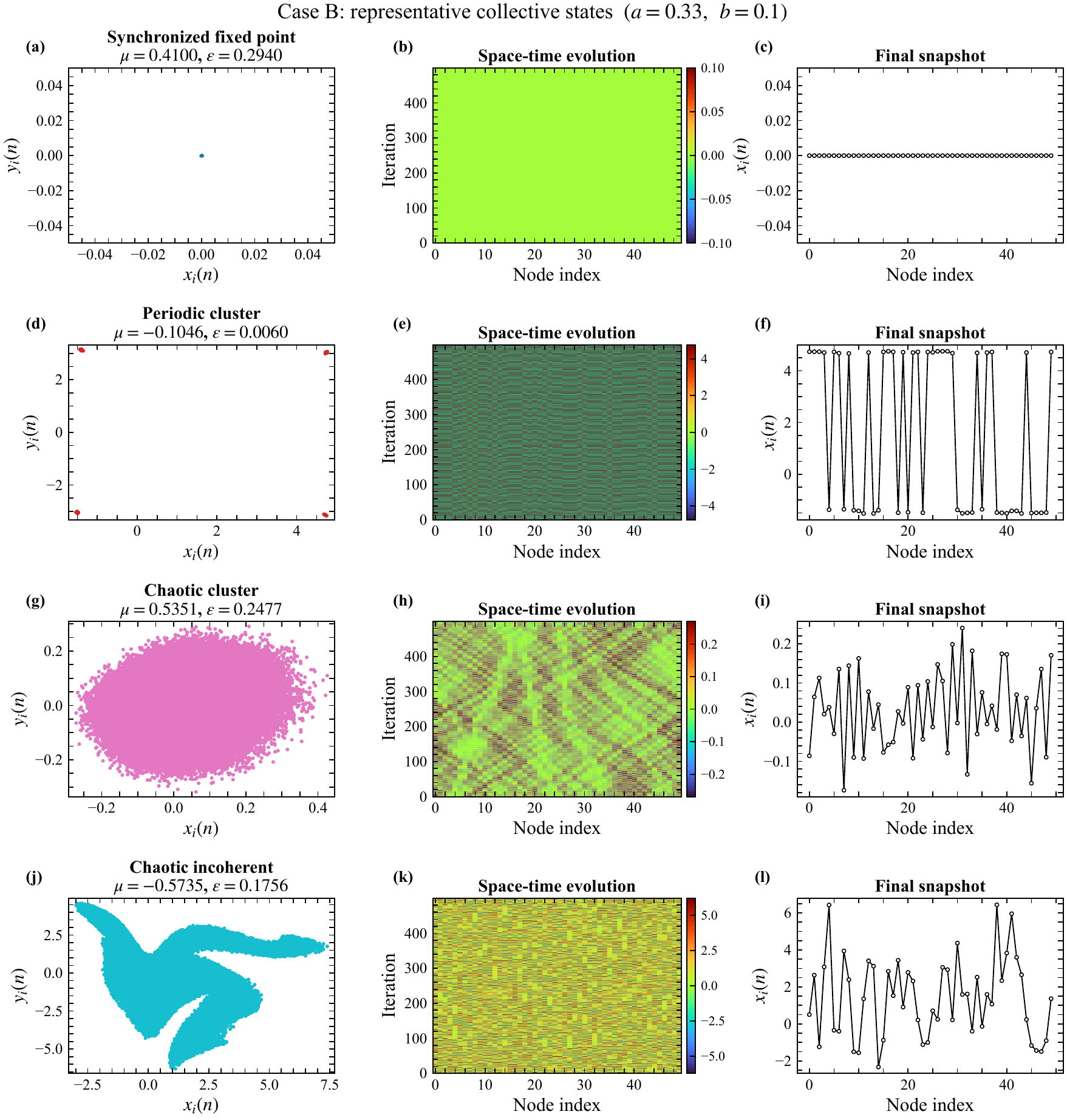} 
\caption{Representative collective states of the coupled Gumowski--Mira network for Case~B, with $a=0.33$ and $b=0.10$. The rows show a synchronized fixed point at $(\mu,\epsilon)=(0.4100,0.2940)$, a periodic cluster state at $(\mu,\epsilon)=(-0.1046,0.0060)$, a chaotic cluster state at $(\mu,\epsilon)=(0.5351,0.2477)$, and a chaotic-incoherent state at $(\mu,\epsilon)=(-0.5735,0.1756)$. The first column shows the physical $(x_i,y_i)$ phase-space trajectories of all nodes, the second column shows the final 500 iterations of the space-time evolution, with color representing the instantaneous value of $x_i(n)$, and the third column shows the final network snapshot.}
\label{fig:representative_states_caseB} 
\end{figure}

The collective-state diagrams identify synchronized-periodic, periodic-cluster, chaotic-cluster, and incoherent regimes. Figures~\ref{fig:representative_states_caseA} and \ref{fig:representative_states_caseB} show representative states for Cases~A and~B, respectively. For each case, the rows correspond to a synchronized fixed point, periodic cluster, chaotic cluster, and chaotic-incoherent state. The physical phase-space trajectories, space-time evolution, and final network configuration provide complementary views of their temporal and spatial organization.

The synchronized-periodic regime is characterized by regular temporal dynamics and strong global coherence. In this state, the synchronization error is small, the network forms a single coherent group, and the largest Lyapunov exponent is negative. For the representative states shown in Figs.~\ref{fig:representative_states_caseA} and \ref{fig:representative_states_caseB}, all nodes converge to the origin. These states are therefore identified more specifically as synchronized fixed points. Their phase-space trajectories and final network snapshots collapse to the origin, while the space-time panels remain spatially uniform.

The periodic cluster-synchronization regime is also temporally regular, but the network does not collapse into a single synchronized group. Instead, the nodes organize into a small number of detected clusters. The selected periodic-cluster representatives have negative finite-time largest Lyapunov exponents, synchronization errors above $E_{\mathrm{tol}}$, and rounded median cluster number $N_c=2$. This state is particularly important because periodic cluster synchronization occupies substantial regions of the parameter plane in both Cases A and B, showing that clustering is a robust emergent feature of the coupled Gumowski--Mira network. The two cases exhibit different instantaneous spatial arrangements, illustrating that the same periodic-cluster classification can contain different cluster geometries.

The chaotic cluster-synchronization regime combines irregular temporal dynamics with spatial organization. In this regime, the finite-time largest Lyapunov exponent is positive, while the number of detected clusters remains small. The selected Case~A and Case~B chaotic-cluster representatives both have rounded median cluster number $N_c=2$. Their finite-time largest Lyapunov exponents are $1.323\times10^{-1}$ and $6.120\times10^{-3}$, respectively. The structured space-time panels show that partial spatial organization can persist even when the phase-space trajectories are irregular. 

The Case~B value is comparatively close to the adopted finite-time threshold $\lambda_{\mathrm{tol}}=10^{-3}$. Consequently, the finite-time exponent alone might suggest only weak chaotic instability. The complete spectrum reported below, however, identifies 21 exponents above $10^{-3}$ and therefore provides stronger evidence of high-dimensional instability. 

Finally, the incoherent regime is characterized by irregular temporal dynamics and weak spatial organization. In this state, the number of detected clusters is large, the phase-space trajectories become broadly distributed, the space-time plot loses low-cluster organization, and the final spatial snapshot displays a fragmented node configuration. 

The selected chaotic-incoherent states have positive finite-time largest Lyapunov exponents and rounded median cluster numbers $N_c=41$ for Case~A and $N_c=36$ for Case~B. These values are substantially larger than the fragmentation threshold $N_c^*=3$ and distinguish the chaotic-incoherent states from the two chaotic-cluster representatives. 

The reported $N_c$ values are rounded medians of cluster counts evaluated at multiple post-transient sampling times. A single final network snapshot therefore need not display exactly the reported number of amplitude groups. This distinction is particularly relevant for chaotic clusters because the instantaneous node amplitudes and apparent cluster memberships can fluctuate in time. 

\subsection{Complete Lyapunov spectra} \label{sec:complete_lyapunov_spectra} 

To distinguish nonchaotic dynamics, ordinary chaos, and hyperchaos, we computed the complete 100-dimensional Lyapunov spectrum for all eight representative states using QR orthonormalization. The results are reported in Table~\ref{tab:hyperchaos_check_appendix}. A negative leading Lyapunov exponent confirms nonchaotic temporal dynamics but does not, by itself, establish complete synchronization. Table~V therefore also reports the median synchronization error \(\widetilde{E}\) and median cluster number \(\widetilde{N}_c\) obtained from the 50-initial-condition verification. Complete synchronization is confirmed independently by \(\widetilde{E}\leq E_{\mathrm{tol}}\) and \(\widetilde{N}_c=1\). The agreement percentage indicates how frequently each requested representative classification was recovered from independent initial network states. The synchronized fixed points and periodic clusters are nonchaotic in both local parameter cases. All 100 exponents lie below the positive threshold of $10^{-3}$ for these four representatives. 

Here, ``chaotic cluster'' denotes the finite-time collective-state category defined by $\lambda_{\max}>\lambda_{\mathrm{tol}}$, $E>E_{\mathrm{tol}}$, and $2\leq N_c\leq N_c^\ast$. By contrast, ``hyperchaotic'' describes the complete-spectrum result $n_+\geq2$ for the same trajectory. Thus, the cluster label characterizes spatial organization, whereas the hyperchaotic label characterizes the number of unstable tangent-space directions.

For Case~A, the synchronized fixed point has leading exponent $\lambda_1=-1.691\times10^{-1}$, while the periodic cluster has $\lambda_1=-3.815\times10^{-3}$. For Case~B, the corresponding values are $-3.147\times10^{-2}$ and $-4.405\times10^{-1}$. These complete-spectrum results are consistent with the negative finite-time largest exponents used in the collective-state classification. 

Both chaotic-cluster representatives are hyperchaotic. The Case~A chaotic cluster at $(\mu,\epsilon)=(-0.5760,0.2928)$ has 
\[\lambda_1=1.592\times10^{-1}, \qquad n_+=19,\]
whereas the Case~B chaotic cluster at $(\mu,\epsilon)=(0.5351,0.2477)$ has \[\lambda_1=2.604\times10^{-3}, \qquad n_+=21.\]
Thus, low-cluster spatial organization coexists with expansion along many independent tangent-space directions in both cases. Although the leading exponent of the Case~B chaotic cluster is comparatively small, its 21 positive exponents show that the hyperchaotic classification is not determined by a single marginally positive direction. 

The chaotic-incoherent representatives exhibit still more broadly distributed instability. The Case~A chaotic-incoherent state has 
\[\lambda_1=3.761\times10^{-1}, \qquad n_+=50,\] while the Case~B chaotic-incoherent state has \[\lambda_1=4.071\times10^{-1}, \qquad n_+=50.\] The 50 positive exponents in each 100-dimensional spectrum show that instability is distributed across half of the tangent-space directions under the adopted positive-exponent threshold.

The complete spectra therefore reveal a hierarchy of collective dynamics. The synchronized fixed points and periodic clusters are nonchaotic. The chaotic clusters retain low-cluster spatial organization but are hyperchaotic, with 19 and 21 positive exponents. The chaotic-incoherent states are also hyperchaotic and have 50 positive exponents in both cases. The distinction between chaotic clustering and chaotic incoherence is therefore primarily spatial, but it is accompanied by a substantial difference in the number of unstable tangent-space directions. 

The signs of the leading exponents and the positive-exponent counts were monitored at intermediate accumulation times between 2000 and 20000 iterations. The leading complete-spectrum exponents were also compared with the corresponding two-trajectory estimates. These checks support the reported nonchaotic and hyperchaotic classifications.

\subsection{Initial-condition dependence of selected states}  \label{sec:network_initial_condition_dependence}
\textcolor{black}{To assess the robustness of the representative collective-state classifications, each of the eight parameter pairs shown in Figs.~\ref{fig:representative_states_caseA} and \ref{fig:representative_states_caseB} was recomputed using 50 independent network initial conditions. For each realization, the initial states of the $N=50$ nodes were sampled independently from $\mathcal{U}[-1,1]$. The finite-time largest Lyapunov exponent, synchronization error, rounded median cluster number, and collective-state classification were evaluated from the same reference trajectory. The analysis therefore comprised 400 network realizations, all of which remained numerically valid.}

\textcolor{black}{Six of the eight representative parameter pairs retained their requested classifications for all sampled initial conditions. The two exceptions occurred in Case A. At the periodic-cluster parameter pair $(\mu,\varepsilon)=(0.4918,0.0174)$, 10\% of the realizations retained the periodic-cluster classification, whereas 90\% converged to synchronized-periodic dynamics. The largest Lyapunov exponent remained negative for all realizations, indicating coexistence between different spatial organizations of nonchaotic dynamics. At the chaotic-cluster parameter pair $(\mu,\varepsilon)=(-0.5760,0.2928)$, 92\% of the realizations retained the chaotic-cluster classification, whereas 8\% converged to synchronized-periodic dynamics. This parameter pair therefore supports both chaotic low-cluster and regular synchronized outcomes under different initial conditions.}

\textcolor{black}{The synchronized fixed-point and chaotic-incoherent parameter pairs in Case A retained their requested classifications for all sampled initial conditions. All four Case B parameter pairs also retained their requested classifications under the same initialization protocol. These results show that the robustness of a collective classification depends on the selected parameter region. In particular, the two Case A cluster parameter pairs exhibit multiple collective outcomes under different initial conditions. The collective-state diagrams should therefore be interpreted as classifications obtained under the stated initialization protocol rather than as exhaustive basin diagrams of the full $2N$-dimensional network state space. The results are summarized in Table~\ref{tab:initial_condition_robustness}.}

\begin{table}[t]
\centering 
\caption{Robustness of the eight representative collective states over 50 independent network initial conditions per parameter pair. Agreement denotes the percentage of realizations that retained the requested finite-time classification.} 
\label{tab:initial_condition_robustness}
\begin{tabular}{llcl} 
\hline 
Case & State & Agreement & Other outcomes \\ 
\hline 
A & Synchronized fixed point & $100\%$ & None \\ 
A & Periodic cluster & $10\%$ & $90\%$ synchronized periodic \\ 
A & Chaotic cluster & $92\%$ & $8\%$ synchronized periodic \\
A & Chaotic incoherent & $100\%$ & None \\ B & Synchronized fixed point & $100\%$ & None \\
B & Periodic cluster & $100\%$ & None \\ B & Chaotic cluster & $100\%$ & None \\ 
B & Chaotic incoherent & $100\%$ & None \\ 
\hline 
\end{tabular}
\end{table}

\subsection{Initial-condition dependence of the isolated map} \label{sec:isolated_initial_condition_dependence}
\textcolor{black}{The initial-condition dependence of the six representative isolated-map states shown in Fig.~\ref{fig:phase_portraits} was examined using 100 independently sampled initial conditions from $[-1,1]^2$. For each initial condition, the finite-time largest Lyapunov exponent was computed using the same transient and averaging protocol adopted in the main analysis. The attracting-periodic and chaotic representatives were robust under this test: all sampled trajectories retained their corresponding regular or chaotic classifications as the Lyapunov averaging interval was increased.}

\begin{figure}[t]
\centering 
\includegraphics[width=\columnwidth]{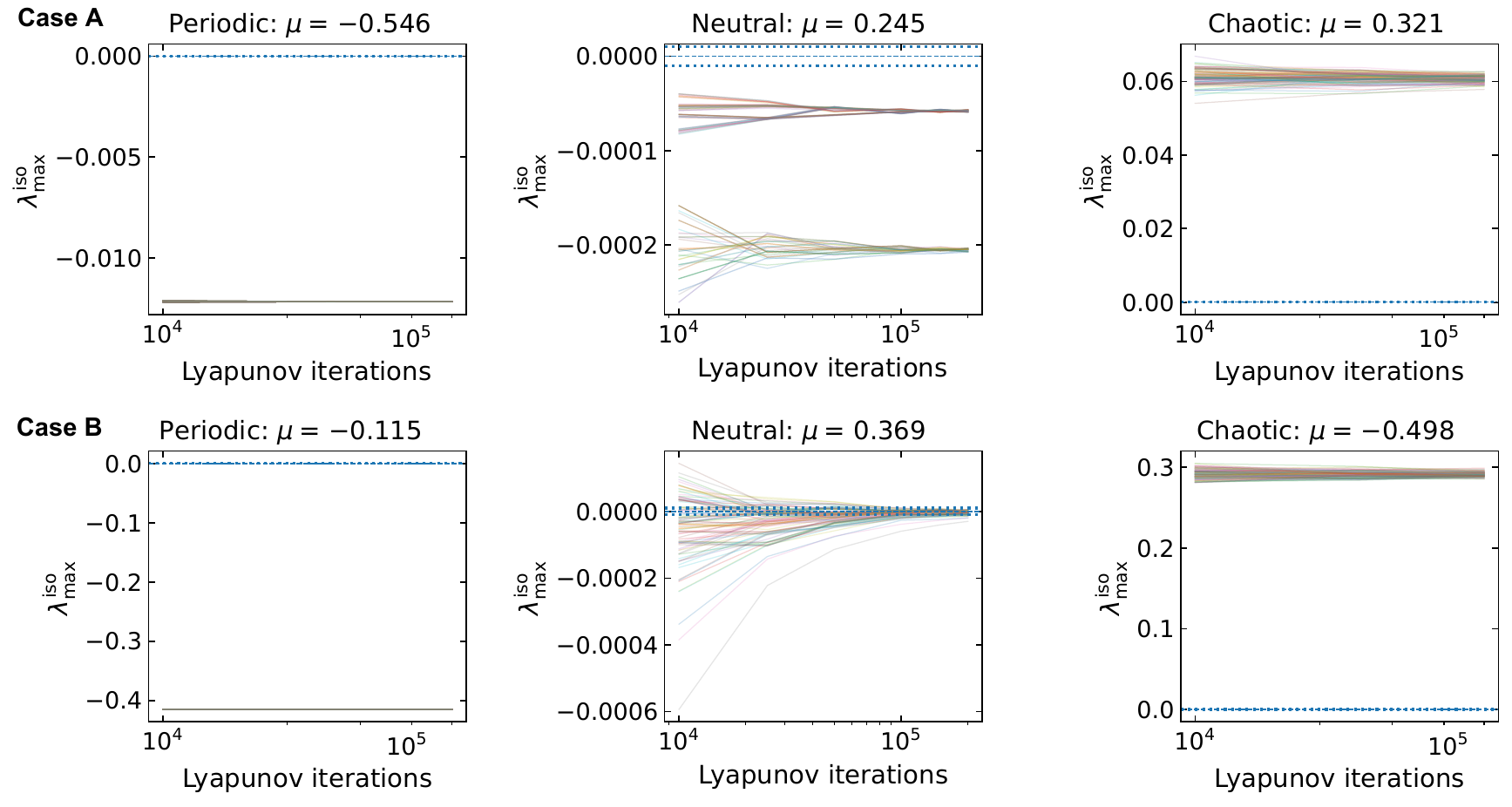}
\caption{Convergence of the finite-time largest Lyapunov exponent with the averaging interval for the six representative isolated-map states, using 100 independently sampled initial conditions for each parameter value. The attracting-periodic and chaotic states retain clearly negative and positive values, respectively, whereas the numerically neutral states show slower convergence toward zero.}
\label{fig:100ICs} 
\end{figure}

\textcolor{black}{The convergence of the largest Lyapunov exponent with the averaging interval is summarized in Fig.~\ref{fig:100ICs}. The attracting-periodic representatives remained clearly negative, whereas the chaotic representatives remained positive over the tested intervals. By contrast, the numerically neutral representatives converged more slowly toward zero and showed greater variability across initial conditions. This behaviour was particularly pronounced for Case B at $\mu=0.3690$, where the finite-time estimates approached zero from both positive and negative values.} 

\begin{figure}[t]
\centering 
\includegraphics[width=0.9\columnwidth]{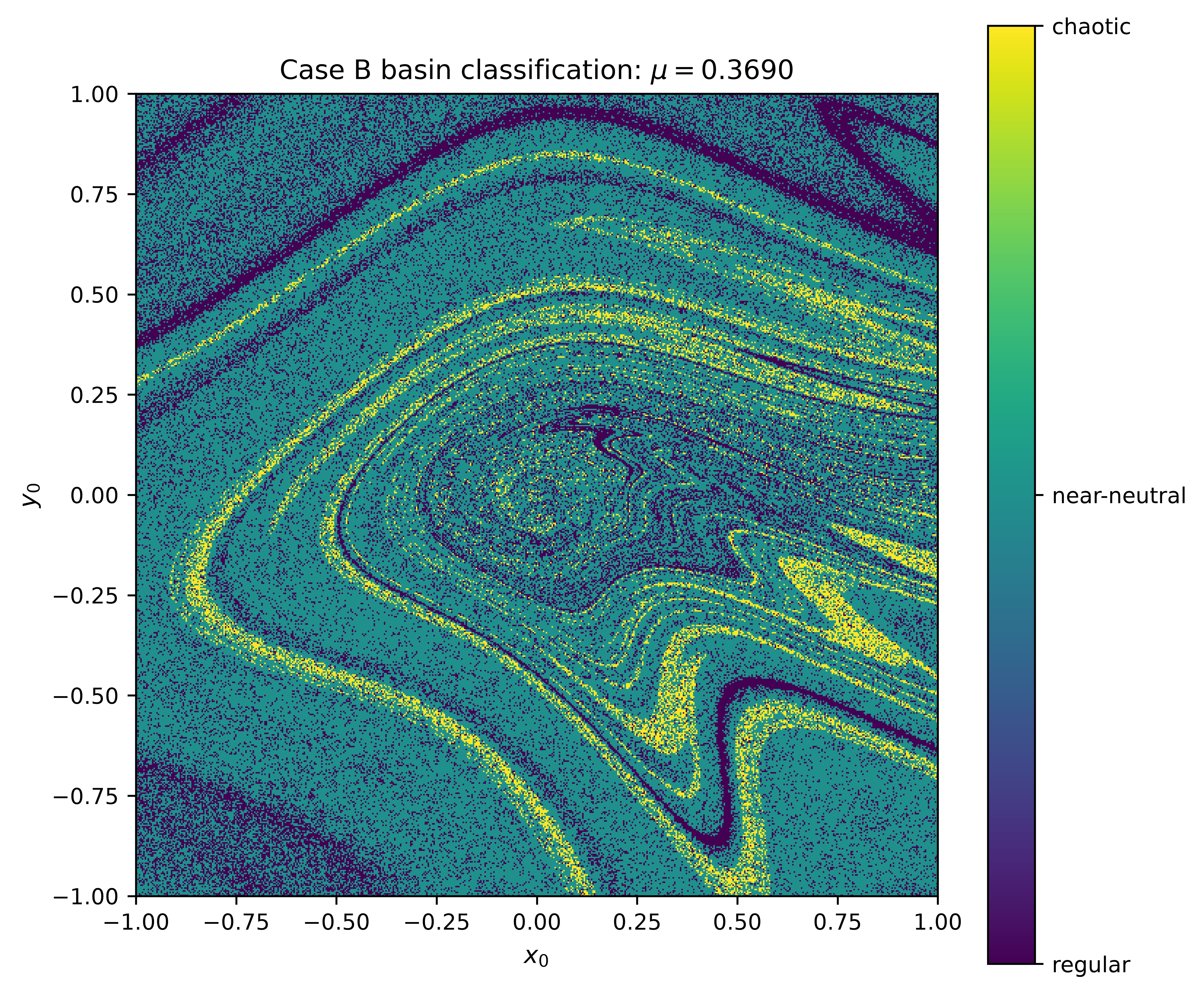} 
\caption{Finite-time initial-condition classification for the isolated Case B map at $\mu=0.3690$, evaluated over a $500\times500$ grid in $[-1,1]^2$. Initial conditions are classified as contracting, near-neutral, or expanding according to the finite-time largest Lyapunov exponent, using $\lambda_{\mathrm{tol}}=10^{-3}$. The structured regions describe differences in finite-time convergence and should not be interpreted as confirmed basins of distinct asymptotic attractors.}
\label{fig:basin_caseB_mu0369} 
\end{figure}

\textcolor{black}{To investigate the apparent initial-condition dependence of the near-neutral Case B representative more closely, the isolated map at $\mu=0.3690$ was evaluated over a $500\times500$ grid of initial conditions. Using the tolerance $\lambda_{\mathrm{tol}}=10^{-3}$, each initial condition was classified as contracting, near-neutral, or expanding according to whether the computed finite-time largest Lyapunov exponent satisfied $\lambda_{\max}<-\lambda_{\mathrm{tol}}$, $|\lambda_{\max}|\leq\lambda_{\mathrm{tol}}$, or $\lambda_{\max}>\lambda_{\mathrm{tol}}$, respectively. As shown in Fig.~\ref{fig:basin_caseB_mu0369}, these finite-time classifications form structured regions in the initial-condition plane.}

\textcolor{black}{Representative trajectories selected from the different regions in Fig.~\ref{fig:basin_caseB_mu0369} approached near-zero largest Lyapunov exponents when longer averaging intervals were used and occupied similar regions of phase space. The apparent variation in Fig.~\ref{fig:basin_caseB_mu0369} is therefore more consistent with differences in finite-time convergence near a neutral regime than with confirmed basins of distinct asymptotic attractors.}

\textcolor{black}{Overall, the attracting-periodic and chaotic isolated-map representatives were robust under the sampled initialization protocol. By contrast, states with largest Lyapunov exponents close to zero require careful interpretation because their finite-time classifications can depend on the initial condition, averaging interval, and numerical tolerance.}

\subsection{Finite-time and tolerance sensitivity} \label{sec:finite_time_sensitivity} 
\textcolor{black}{Because the collective-state classifications are based on finite-time numerical diagnostics, we examined their sensitivity to the post-transient observation interval and the numerical thresholds used in the classification. The observation interval was varied over
\[ T_{\mathrm{obs}}\in\{1000,2000,5000,10000,20000\} \] 
iterations, following the same transient-removal procedure. We also tested the discrete threshold values 
\[ \lambda_{\mathrm{tol}} \in\{10^{-4},10^{-3},10^{-2}\}, \qquad E_{\mathrm{tol}} \in\{10^{-10},10^{-8},10^{-6}\}, \] 
and \[ \eta\in\{0.025,0.05,0.075\}. \] }

\textcolor{black}{The synchronized fixed-point and chaotic-incoherent representatives retained their classifications over all tested observation intervals. The Case A chaotic-cluster and Case B periodic-cluster representatives also retained their median-based classifications. By contrast, the Case A parameter pair originally selected as a periodic-cluster state was classified as synchronized periodic at every tested observation interval. This result agrees with the multiple-initial-condition analysis in Sec.~\ref{sec:representative_collective_states}.B, in which synchronized-periodic dynamics were the dominant outcome. The parameter pair is therefore interpreted as multistable rather than as a uniformly robust periodic-cluster state.}

\textcolor{black}{The Case B chaotic-cluster parameter pair exhibited greater finite-time structural sensitivity. Its largest Lyapunov exponent remained positive over all tested observation intervals, confirming the persistence of chaotic temporal dynamics. However, its rounded median cluster number satisfied $N_c\geq 2$ for observation intervals up to 5000 iterations and decreased to $N_c=1$ for intervals of 10000 and 20000 iterations. Because the synchronization error remained above $E_{\mathrm{tol}}$, the trajectory did not satisfy the criterion for complete synchronization at the two longest observation intervals. It also no longer satisfied the adopted chaotic-cluster criterion, which requires $2\leq N_c\leq N_c^*$. Thus, the temporal chaotic character persisted, whereas the detected spatial organization depended on the observation interval. }

\textcolor{black}{The synchronization tolerance did not alter the representative classifications because the synchronization errors were well separated from the tested values of $E_{\mathrm{tol}}$. The Lyapunov tolerance affected only states with finite-time largest Lyapunov exponents close to zero. In particular, a weakly contracting synchronized trajectory with $\lambda_{\max}\approx-3.82\times10^{-3}$ was classified as synchronized periodic for $\lambda_{\mathrm{tol}}=10^{-4}$ and $10^{-3}$, but as synchronized neutral for $\lambda_{\mathrm{tol}}=10^{-2}$. Similarly, a weakly expanding trajectory with $\lambda_{\max}\approx2.66\times10^{-3}$ was classified as temporally neutral when the largest Lyapunov tolerance was used. }

 \textcolor{black}{The cluster threshold primarily affected the number of detected clusters in chaotic and strongly fragmented states because $\eta$ controls the separation of adjacent sorted node amplitudes. We therefore retained \[ \lambda_{\mathrm{tol}}=10^{-3}, \qquad E_{\mathrm{tol}}=10^{-8}, \qquad \eta=0.05 \] as the baseline values. Overall, the principal classifications were stable under the tested numerical settings, while classifications close to temporal or structural boundaries showed sensitivity to the observation interval or numerical threshold. Such cases were interpreted using the underlying values of $\lambda_{\max}$, $E$, and $N_c$, rather than from the categorical label alone.}

 \subsection{Network-size dependence} \label{sec:network_size} 
 
 \textcolor{black}{The influence of network size was examined using $N\in\{25,50,100\}$, with ten independently generated initial network states for each representative parameter pair and network size. All tested realizations remained numerically valid.}
 
 \textcolor{black}{The synchronized fixed-point representatives retained their classifications for all sampled initial conditions at the three network sizes. The chaotic-incoherent representatives were also qualitatively robust, although their absolute rounded median cluster numbers increased with $N$. For Case A, the median cluster numbers were $23$, $41$, and $66.5$ for $N=25$, $50$, and $100$, respectively; the corresponding values for Case B were $21$, $36$, and $52$. However, the normalized cluster number $N_c/N$ decreased as $N$ increased. Thus, the absolute cluster count is not a size-independent measure, even when the qualitative incoherent classification is retained. }
 
 \textcolor{black}{The Case B periodic-cluster representative retained regular temporal dynamics and low-cluster organization at all three network sizes. By contrast, the Case A parameter pair originally selected as a periodic-cluster representative was predominantly synchronized periodic. The periodic-cluster outcome occurred in none of the ten realizations for $N=25$ and in one of the ten realizations for each of $N=50$ and $N=100$. This behaviour is consistent with the initial-condition analysis in Sec.~\ref{sec:representative_collective_states}.B and supports the interpretation of this parameter pair as multistable rather than as a network-size-robust periodic-cluster state.}
 
 \textcolor{black}{The chaotic-cluster representatives exhibited greater finite-size sensitivity. For Case A, the chaotic-cluster classification occurred in $30\%$, $100\%$, and $90\%$ of the realizations for $N=25$, $50$, and $100$, respectively. For Case B, the corresponding fractions were $100\%$, $90\%$, and $10\%$. At $N=100$, the Case B trajectories generally retained positive largest Lyapunov exponents, but most had $N_c=1$ while their synchronization errors remained above $E_{\mathrm{tol}}$. These trajectories therefore retained chaotic temporal dynamics but no longer satisfied either the adopted chaotic-cluster criterion or the criterion for complete synchronization.}
 
 \textcolor{black}{Overall, the synchronized fixed-point, chaotic-incoherent, and Case B periodic-cluster representatives were qualitatively robust over the tested network sizes. In contrast, the prevalence and detected spatial organization of the selected chaotic-cluster states depended on $N$. We therefore interpret $N=50$ as a representative baseline network size and do not claim size-independent cluster counts or universal finite-size robustness of the selected clustered states.}

\section{Conclusion}\label{sec:VIII}

In this work, we investigated the collective dynamics of \myEdit{a nearest-neighbor ring of} coupled Gumowski--Mira maps under two intrinsic local parameter regimes. The largest Lyapunov exponent, synchronization error, \myEdit{rounded median cluster number}, and combined collective-state classification reveal synchronized, clustered, and incoherent regimes whose locations and distributions depend on the local parameters $(a,b)$. The collective parameter spaces show that increasing coupling can promote attracting regular behavior and synchronization, although the resulting organization differs between the two local regimes.

Case~A exhibits \myEdit{an overall sloping but nonuniform transition toward regular and coherent collective dynamics as the coupling strength increases.} Case~B displays a more heterogeneous \myEdit{organization, with broad incoherent regions, localized clustered structures, and an extended synchronized-periodic domain at strong coupling and positive $\mu$.} 

These differences demonstrate that changing the intrinsic map dynamics modifies the organization of the collective parameter space even when the network topology and coupling rule are unchanged. 

\myEdit{Long-run representative simulations confirm synchronized fixed points, periodic two-cluster states, chaotic two-cluster states, and chaotic-incoherent dynamics in both cases. The synchronized representatives converge to the origin and have negative leading Lyapunov exponents, numerically vanishing synchronization errors, and rounded median cluster number $N_c=1$. The periodic-cluster representatives have negative leading exponents and rounded median cluster number $N_c=2$. The chaotic-cluster representatives combine positive finite-time largest Lyapunov exponents with low-cluster spatial organization, whereas the chaotic-incoherent representatives combine positive exponents with substantially larger rounded median cluster numbers.} 

\myEdit{ The complete Lyapunov spectra provide a clear dynamical distinction among these states. The synchronized fixed points and periodic clusters are nonchaotic, with no Lyapunov exponents above the adopted positive threshold of $10^{-3}$. By contrast, the Case~A and Case~B chaotic-cluster representatives are hyperchaotic, with 19 and 21 positive exponents, respectively. The chaotic-incoherent representatives are also hyperchaotic and have 50 positive exponents in both cases. These results demonstrate that low-cluster spatial organization can coexist with instability along many independent tangent-space directions. The larger positive-exponent counts of the chaotic-incoherent states indicate more broadly distributed instability than in the chaotic clusters.} 

\textcolor{black}{The multiple-initial-condition analysis showed that six of the eight representative parameter pairs retained their requested classifications for all sampled initial conditions. The two exceptions occurred in Case A, where the selected periodic- and chaotic-cluster parameter pairs exhibited multistability involving synchronized-periodic outcomes. Consequently, the collective parameter diagrams should be interpreted as realization-based classifications under the stated initialization protocol rather than as exhaustive basin diagrams. The combined use of finite-time phase classification, long-run representative trajectories, complete Lyapunov spectra, and multiple-initial-condition verification provides a more complete description of the network dynamics than any one diagnostic alone.} 

\textcolor{black}{Future work may examine basin structure, parameter mismatch, noise, and the temporal evolution of cluster membership in hyperchaotic clustered states. It would also be valuable to determine how the number of positive Lyapunov exponents scales with network size and to extend the comparative analysis of collective dynamics across the $(a,b)$ plane for different fixed values of $(\mu,\varepsilon)$. Finally, since the present conclusions are specific to nearest-neighbor ring coupling, future studies should investigate how nonlocal, global, and heterogeneous coupling architectures modify synchronization thresholds, cluster organization, basin prevalence, and the boundaries between collective states.}


\appendix

\section{Numerical Diagnostics and Supplementary Tables}
\label{app:appendix}

This appendix provides additional details on the numerical diagnostics used to characterize the collective dynamics of the coupled Gumowski--Mira network. The diagnostics include the largest Lyapunov exponent, synchronization error, cluster-count measure, and collective-state classification used to construct the phase diagrams. Supplementary tables summarize the local parameter regimes, cluster-count interpretation, phase-classification criteria, and \myEdit{complete Lyapunov-spectrum calculation}.

\subsection{Largest Lyapunov exponent}

The largest Lyapunov exponent is used to distinguish regular collective dynamics from irregular or chaotic collective dynamics. For a discrete-time system, the Lyapunov exponent measures the average exponential rate at which two nearby trajectories separate. If the largest Lyapunov exponent is negative, nearby trajectories converge and the asymptotic dynamics are regular. If it is positive, the system displays sensitive dependence on initial conditions.

For the coupled network, we compute the largest Lyapunov exponent numerically by evolving two initially close network states. Let 
\begin{equation}
\myEdit{ \mathbf{Z}(n)= \left[ x_1(n),\ldots,x_N(n), y_1(n),\ldots,y_N(n) \right]^{\mathsf T} } 
\label{eq:network_state_vector}
\end{equation}
denote the full $2N$-dimensional network state at iteration $n$. A perturbed trajectory 
\begin{equation} 
\mathbf{Z}'(0) = \mathbf{Z}(0)+\delta_0\mathbf{v} 
\end{equation}
is initialized with $\delta_0\ll1$, where $\mathbf{v}$ is a random unit perturbation vector. After each iteration, the separation between the reference and perturbed trajectories is
\begin{equation} 
d(n) = \left\| \mathbf{Z}'(n)-\mathbf{Z}(n) \right\|. 
\end{equation}
The perturbation is then rescaled to length $\delta_0$ while preserving its direction. \myEdit{ After discarding $n_0$ transient iterations, the finite-time largest Lyapunov exponent is estimated over $T_\lambda$ renormalization steps as
\begin{equation} 
\lambda_{\max} = \frac{1}{T_\lambda} \sum_{k=0}^{T_\lambda-1} \ln\left( \frac{d(n_0+k+1)}{\delta_0} \right). 
\label{eq:lyapunov_appendix} 
\end{equation} }

Throughout this work, parameter values satisfying \begin{equation} \lambda_{\max}<-\lambda_{\mathrm{tol}} \end{equation} are interpreted as regular collective dynamics, while values satisfying \begin{equation} \lambda_{\max}>\lambda_{\mathrm{tol}} \end{equation} are interpreted as irregular or chaotic collective dynamics according to the finite-time largest Lyapunov exponent. Values satisfying \begin{equation} |\lambda_{\max}| \leq \lambda_{\mathrm{tol}} \end{equation} are treated as numerically neutral or marginal. \myEdit{ The finite-time largest Lyapunov exponent is used in the parameter-space classification because computing the complete $2N$-dimensional spectrum at every parameter pair would be substantially more expensive. Complete spectra are therefore calculated only for the eight representative states. }

\subsection{Synchronization error} 
To quantify global synchronization, we compute the dispersion of the node $x$-variables about their network mean. At iteration $n$, the network mean is
\begin{equation} 
\bar{x}(n) = \frac{1}{N} \sum_{i=1}^{N}x_i(n). 
\label{eq:network_mean} 
\end{equation} 

The instantaneous synchronization error is 
\begin{equation} E(n) = \frac{1}{N} \sum_{i=1}^{N} \left[ x_i(n)-\bar{x}(n) \right]^2. 
\label{eq:instantaneous_sync_error} 
\end{equation} 

After discarding the first $n_0=1500$ transient iterations, the time-averaged synchronization error is evaluated over the subsequent $T=2000$ iterations as 
\begin{equation} 
E = \frac{1}{T} \sum_{n=n_0}^{n_0+T-1} E(n). 
\label{eq:time_averaged_sync_error} 
\end{equation} 

Complete synchronization is identified when 
\begin{equation} E\leq E_{\mathrm{tol}}, \qquad E_{\mathrm{tol}}=10^{-8}. 
\label{eq:synchronization_condition} 
\end{equation}

Values of $E$ above this tolerance indicate global desynchronization or partial synchronization. However, a nonzero global synchronization error does not distinguish clustered organization from incoherence. We therefore complement $E$ with the cluster-count diagnostic described below.

\subsection{Cluster-count measure} 
To identify cluster-synchronized states, we calculate the number of distinct groups formed by the node variables after removal of the transient. Let 
\begin{equation} \mathbf{x}(n) = \left[ x_1(n),x_2(n),\ldots,x_N(n) \right]
\end{equation} 
denote the vector of node states at a fixed post-transient iteration \(n\). At each sampled iteration, the node values are sorted in ascending order,
\begin{equation} x_{(1)}(n) \leq x_{(2)}(n) \leq \cdots \leq x_{(N)}(n), 
\end{equation} 
where \(x_{(j)}(n)\) denotes the \(j\)-th ordered node value. A boundary between two clusters is introduced whenever the gap between adjacent sorted values exceeds the prescribed tolerance \(\eta\). 
The instantaneous cluster number is therefore defined by 
\begin{equation} 
N_c(n) = 1 + \sum_{j=1}^{N-1} \mathbb{I} \left[ x_{(j+1)}(n)-x_{(j)}(n)>\eta \right], 
\label{eq:cluster_count} 
\end{equation} 
where \(\mathbb{I}[\cdot]\) is the indicator function. In the numerical calculations, we use the fixed threshold \begin{equation} \eta=0.05. \end{equation} \myEdit{ No additional minimum cluster size is imposed by this threshold procedure. Consequently, an isolated node may be retained as a singleton cluster in the raw instantaneous cluster count. However, a state is assigned a cluster-synchronization classification only when the reported cluster number satisfies 
\begin{equation}
2\leq N_c\leq N_c^\ast, 
\end{equation} where \(N_c^\ast=3\) is the fragmentation threshold. Thus, a state with \(N_c=1\) cannot be classified as periodic, neutral, or chaotic cluster synchronization. } \myEdit{ The instantaneous cluster number \(N_c(n)\) is evaluated at 20 uniformly distributed post-transient sampling times. The reported value \(N_c\) is the rounded median of these 20 instantaneous cluster counts. This temporal summary reduces sensitivity to isolated fluctuations in the detected number of clusters. } 

\myEdit{ Because \(N_c\) is a rounded temporal median, the number of distinct amplitude groups visible in a single final network snapshot need not equal the reported value. This distinction is particularly important for chaotic clustered states, for which the instantaneous node amplitudes and apparent cluster memberships may fluctuate over time. The interpretation of \(N_c\) is summarized in Table~\ref{tab:cluster_interpretation_appendix}. } 

\begin{table}[t] 
\centering 
\caption{ \myEdit{ Interpretation of the rounded median cluster number \(N_c\) under the fragmentation threshold \(N_c^\ast=3\) used in the collective-state classification. } } 
\label{tab:cluster_interpretation_appendix} 
\begin{tabular}{cl}
\hline 
Condition & Spatial organization \\ 
\hline \(N_c=1\) & One detected cluster \\ 
\(2\leq N_c\leq N_c^\ast\) & Low-cluster organization \\ 
\(N_c>N_c^\ast\) & Incoherent or highly fragmented organization \\ 
\hline
\end{tabular}
\end{table}


The cluster number alone also does not distinguish regular from chaotic temporal dynamics. A low-cluster state is identified as periodic, neutral, or chaotic by combining \(N_c\) with the global synchronization error \(E\) and the finite-time largest Lyapunov exponent \(\lambda_{\max}\). States satisfying \(N_c>N_c^\ast\) are assigned to the common incoherent category in the combined phase diagram, irrespective of whether their temporal dynamics are regular, numerically neutral, or chaotic. 

Throughout this work, the term ``cluster synchronization'' is retained as an operational finite-time classification describing coherent organization into multiple detected clusters under the stated threshold criterion.




\subsection{Collective-state classification} \label{app:collective_state_classification} 
The largest Lyapunov exponent detects whether the collective dynamics are regular, numerically neutral, or chaotic. The synchronization error quantifies global coherence, and the rounded median cluster number estimates the number of coherent groups. Combining these indicators, we classify the collective dynamics into synchronized-periodic, synchronized-neutral, synchronized-chaotic, periodic-cluster, neutral-cluster, chaotic-cluster, incoherent, and undefined states.

The classification rule is summarized in Table~\ref{tab:collective_classification}. The tolerance $\lambda_{\mathrm{tol}}$ prevents marginal finite-time Lyapunov exponents from being overinterpreted as either regular or chaotic. The tolerance $E_{\mathrm{tol}}$ identifies complete synchronization. The threshold $N_c^\ast=3$ separates low-cluster states from highly fragmented states. 

\myEdit{ The terms ``periodic cluster'' and ``chaotic cluster'' describe different temporal dynamics within the same low-cluster spatial category. A periodic cluster has a negative finite-time largest Lyapunov exponent, whereas a chaotic cluster has a positive finite-time largest Lyapunov exponent. The complete Lyapunov spectrum is subsequently used to determine whether a selected positive-exponent state is ordinarily chaotic or hyperchaotic. } 

\begin{table}[t] \centering \caption{Finite-time collective-state classification.}
\label{tab:collective_classification} 
\begin{tabular}{ll} 
\hline Condition & Collective state \\ 
\hline 
$\lambda_{\max}<-\lambda_{\mathrm{tol}}$, $E\leq E_{\mathrm{tol}}$, $N_c=1$ & Synchronized periodic \\ 
$|\lambda_{\max}|\leq\lambda_{\mathrm{tol}}$, $E\leq E_{\mathrm{tol}}$, $N_c=1$ & Synchronized neutral \\ 
$\lambda_{\max}>\lambda_{\mathrm{tol}}$, $E\leq E_{\mathrm{tol}}$, $N_c=1$ & Synchronized chaotic \\ 
$\lambda_{\max}<-\lambda_{\mathrm{tol}}$, $E>E_{\mathrm{tol}}$, $2\leq N_c\leq N_c^\ast$ & Periodic cluster \\
$|\lambda_{\max}|\leq\lambda_{\mathrm{tol}}$, $E>E_{\mathrm{tol}}$, $2\leq N_c\leq N_c^\ast$ & Neutral cluster \\
$\lambda_{\max}>\lambda_{\mathrm{tol}}$, $E>E_{\mathrm{tol}}$, $2\leq N_c\leq N_c^\ast$ & Chaotic cluster \\ $N_c>N_c^\ast$ & Incoherent \\ 
Invalid or divergent trajectory & Undefined \\ 
\hline 
\end{tabular} 
\end{table} 

The synchronized-periodic category includes fixed points as period-one solutions. The representative synchronized states in Figs.~\ref{fig:representative_states_caseA} and \ref{fig:representative_states_caseB} converge specifically to the origin and are therefore identified as synchronized fixed points.

\myEdit{ The general incoherent category is determined by $N_c>N_c^\ast$, irrespective of the sign of the finite-time largest Lyapunov exponent. An incoherent state is described more specifically as chaotic incoherent when it additionally satisfies $\lambda_{\max}>\lambda_{\mathrm{tol}}$. Both incoherent representatives selected in this work satisfy this additional condition. }

The classification in Table~\ref{tab:collective_classification} is a finite-time numerical classification based on the largest Lyapunov exponent, synchronization error, and rounded median cluster number. Complete synchronization is identified by vanishing synchronization error and one cluster. \myEdit{In the operational classification used here, cluster synchronization
corresponds to finite-time organization into multiple detected clusters, with
\[
E>E_{\mathrm{tol}},
\qquad
2\leq N_c\leq N_c^*.
\]
} 

\myEdit{ Because a positive finite-time largest Lyapunov exponent may be sensitive to the averaging interval, particularly near a transition boundary, the complete Lyapunov spectrum is computed for all eight representative states. This calculation provides a longer-run assessment of the number of unstable tangent-space directions. }

\subsection{Complete Lyapunov-spectrum calculation} \label{app:complete_lyapunov_spectrum} 

For each representative state, we compute the complete $2N$-dimensional Lyapunov spectrum using the standard QR-orthonormalization method ~\cite{Benettin1980a,Benettin1980b}. The tangent-space ordering is 
\begin{equation} \delta\mathbf{Z}(n) = \left[ \delta x_1(n),\ldots,\delta x_N(n), \delta y_1(n),\ldots,\delta y_N(n) \right]^{\mathsf T}. 
\label{eq:tangent_state} 
\end{equation}

For compactness, define 
\begin{equation} \mathbf{A}_n = \frac{\partial\mathbf{X}_{n+1}} {\partial\mathbf{X}_n}, \qquad \mathbf{B}_n = \frac{\partial\mathbf{X}_{n+1}} {\partial\mathbf{Y}_n}, \label{eq:jacobian_blocks} 
\end{equation} 
and 
\begin{equation} 
\mathbf{C}_n = \operatorname{diag} \left[ F_x(x_1(n+1),\mu),\ldots, F_x(x_N(n+1),\mu) \right]. 
\label{eq:C_matrix} 
\end{equation}

The one-step network Jacobian is 
\begin{equation} 
\mathbf{J}_n = 
\begin{pmatrix} \mathbf{A}_n & \mathbf{B}_n \\ -\mathbf{I}+\mathbf{C}_n\mathbf{A}_n & \mathbf{C}_n\mathbf{B}_n
\end{pmatrix}. 
\label{eq:network_jacobian} 
\end{equation} 

The nonzero entries of $\mathbf{A}_n$ are 
\begin{equation} 
(\mathbf{A}_n)_{ii} = (1-\epsilon)F_x(x_i(n),\mu), \qquad (\mathbf{A}_n)_{i,i\pm1} = \frac{\epsilon}{2}, 
\label{eq:A_entries} 
\end{equation} 
where the indices are interpreted periodically. The matrix $\mathbf{B}_n$ is diagonal, with 
\begin{equation} (\mathbf{B}_n)_{ii} = (1-\epsilon) \left[ 1+a\left(1-3b y_i^2(n)\right) \right].
\label{eq:B_entries}
\end{equation} 

Starting from an orthonormal tangent basis, the basis is propagated under $\mathbf{J}_n$ and orthonormalized after every iteration. The Lyapunov exponents are obtained from the accumulated logarithms of the diagonal entries of the resulting triangular factors and are sorted in descending order, 
\begin{equation} \lambda_1\geq\lambda_2\geq\cdots\geq\lambda_{2N}. 
\label{eq:spectrum_ordering} 
\end{equation} 

For $N=50$, each spectrum contains 100 exponents. The reference trajectory is evolved for 5000 transient iterations, and the spectrum is accumulated over $T_{\mathrm{spec}}=20000$ subsequent iterations. Intermediate estimates at 2000, 5000, 10000, 15000, and 20000 iterations are used to monitor convergence. 

Let $n_+$ denote the number of exponents satisfying $\lambda_i>10^{-3}$. A representative state is classified as nonchaotic when $n_+=0$, chaotic when $n_+=1$, and hyperchaotic when $n_+ \geq  2$. This complete-spectrum classification characterizes temporal instability and complements, rather than replaces, the finite-time classification based on $\lambda_{\max}$, $E$, and $N_c$. The results are summarized in Table~\ref{tab:hyperchaos_check_appendix}.

\begin{table*}[!t] 
\centering 
\caption{ Complete Lyapunov-spectrum and collective-state diagnostics for the eight representative parameter pairs shown in Figs.~\ref{fig:representative_states_caseA} and \ref{fig:representative_states_caseB}. Each spectrum contains \(2N=100\) exponents and was accumulated over \(20000\) post-transient iterations. Here, \(\lambda_1\) and \(n_+\) refer to the displayed representative trajectory, where \(n_+\) is the number of exponents satisfying \(\lambda_i>10^{-3}\). The quantities \(\widetilde{E}\) and \(\widetilde{N}_c\) are the median synchronization error and median cluster number, respectively, over 50 independent initial conditions. Agreement is the percentage of valid realizations retaining the requested finite-time classification. Complete synchronization is independently confirmed by \(\widetilde{E}\leq E_{\mathrm{tol}}\) and \(\widetilde{N}_c=1\). } \label{tab:hyperchaos_check_appendix} 
\begin{tabular}{llrrrrrcl} 
\hline 
Case & State & \(\mu\) & \(\epsilon\) & \(\lambda_1\) & \(\widetilde{E}\) & \(\widetilde{N}_c\) & \(n_+\) & Agreement \\ \hline A & Synchronized fixed point & 0.5760 & 0.2940 & \(-1.691\times10^{-1}\) & \(0\) & 1 & 0 & \(100\%\) \\ A & Periodic cluster & 0.4918 & 0.0174 & \(-3.815\times10^{-3}\) & \(8.63\times10^{-19}\) & 1 & 0 & \(10\%\) \\ A & Chaotic cluster & \(-0.5760\) & 0.2928 & \(1.592\times10^{-1}\) & \(2.78\times10^{-2}\) & 2 & 19 & \(92\%\) \\ A & Chaotic incoherent & \(-0.4942\) & 0.0992 & \(3.761\times10^{-1}\) & \(1.07\times10^{1}\) & 41 & 50 & \(100\%\) \\ B & Synchronized fixed point & 0.4100 & 0.2940 & \(-3.147\times10^{-2}\) & \(4.85\times10^{-14}\) & 1 & 0 & \(100\%\) \\ B & Periodic cluster & \(-0.1046\) & 0.0060 & \(-4.405\times10^{-1}\) & \(8.53\) & 3 & 0 & \(100\%\) \\ B & Chaotic cluster & 0.5351 & 0.2477 & \(2.604\times10^{-3}\) & \(1.20\times10^{-2}\) & 2 & 21 & \(100\%\) \\ B & Chaotic incoherent & \(-0.5735\) & 0.1756 & \(4.071\times10^{-1}\) & \(4.43\) & 36 & 50 & \(100\%\) \\ \hline \end{tabular} \vspace{0.5em} \begin{minipage}{0.97\textwidth} \footnotesize The Case~A periodic-cluster parameter pair is multistable. Although the displayed trajectory and complete spectrum correspond to a periodic-cluster realization, only \(10\%\) of the independently initialized realizations retained that classification; the remaining \(90\%\) converged to synchronized-periodic dynamics. Consequently, the median diagnostics for this parameter pair satisfy \(\widetilde{E}\leq E_{\mathrm{tol}}\) and \(\widetilde{N}_c=1\). \end{minipage} \end{table*}

\myEdit{Table~\ref{tab:hyperchaos_check_appendix} shows that the displayed synchronized and periodic-cluster trajectories are nonchaotic. The additional median diagnostics confirm complete synchronization for the synchronized representatives and reveal that the Case~A periodic-cluster parameter pair is multistable, with synchronized-periodic dynamics dominating across independent initial conditions. By contrast, the displayed chaotic-cluster trajectories possess 19 and 21 positive Lyapunov exponents, confirming that their irregular temporal dynamics are hyperchaotic rather than ordinary chaotic. The chaotic-incoherent representatives each possess 50 positive exponents and therefore exhibit instability along a larger number of tangent-space directions.

The displayed Case~B chaotic-cluster trajectory has the smallest positive leading exponent among the four hyperchaotic representatives, \(\lambda_1=2.604\times10^{-3}\). Nevertheless, its 21 exponents above the adopted threshold demonstrate that the hyperchaotic classification does not arise from a single marginally positive tangent-space direction. The supplementary sensitivity analysis shows, however, that its structural cluster classification depends on the initial condition, observation window, and network size. 

These results demonstrate that spatial cluster organization and temporal instability are distinct properties. A network trajectory may exhibit high-dimensional hyperchaotic dynamics while maintaining low-cluster organization over a finite observation interval. Conversely, the temporal hyperchaotic classification may persist even when the detected cluster organization changes under longer observation windows or different network sizes.}

\FloatBarrier

\section*{Acknowledgements} 
\textcolor{black}{\textbf{Hammed O. Fatoyinbo} sincerely acknowledges the support received through the AUT School of Engineering, Computer and Mathematical Sciences 2025 Early Career Research Award.}

\section*{Author Declarations}
\subsection*{Author Contributions}

\textbf{Hammed O. Fatoyinbo}: Conceptualization; Data curation; Formal analysis; Investigation; Methodology; Project administration; Software; Supervision; Validation; Visualization; Writing -- original draft; Writing -- review and editing. \textbf{Indranil Ghosh}: Conceptualization; Formal analysis; Investigation; Methodology; Project administration; Validation; Writing -- original draft; Writing -- review and editing.

\subsection*{Conflict of Interest}
The authors have no conflicts to disclose.


\section*{Code Availability}
\textcolor{black}{All analyses reported in this study were conducted in Python using Jupyter Notebook. The corresponding source code is publicly accessible at \href{https://github.com/hamfat/Gumowski-Mira-Map}{Gumowski-Mira-Map}.}


\end{document}